\documentclass[a4paper,11pt]{article}

\usepackage{amssymb,amsmath,diagrams,times,theorem}
\usepackage{pstcol,pst-node,pst-grad}

\newcommand{\C}{\mathbb{C}}

\newcommand{\N}{\mathbb{N}}

\newcommand{\Z}{\mathbb{Z}}

\newcommand{\Oscr}{{\mathcal O}}

\newcommand{\wis}[1]{{\text{\em \usefont{OT1}{cmtt}{m}{n} #1}}}
 \newcommand{\proofhead}[1]{\par\pagebreak[1]\noindent{\bf#1.\ }}
 \newcommand{\pf}{\proofhead{Proof}}
  \newcommand{\qed}{{\unskip\nolinebreak[1]\hspace{1.5em}\mbox{}\nolinebreak
    \hfill$\Box$\parfillskip=0pt\finalhyphendemerits=0\par\pagebreak[1]}}

\newenvironment{proof}{\pf}{\qed}

\newtheorem{theorem}{Theorem}
\newtheorem{proposition}{Proposition}

{\theorembodyfont{\rmfamily} 
\newtheorem{example}{Example}

\newtheorem{definition}{Definition} }

\title{Non-commutative covers 
and the modular group}

\author{Lieven Le Bruyn and Jan Adriaenssens \\
Departement Wiskunde en Informatica \\ Universiteit Antwerpen  \\
B-2020 Antwerp (Belgium) \\
lieven.lebruyn@ua.ac.be and  jan.adriaenssens@ua.ac.be}

\date{}

\begin{document}

\maketitle

\begin{abstract}
We use non-commutative geometry to study the bulk of finite dimensional representations of the modular group $\Gamma = SL_2(\Z)$. We give specific $2n$-dimensional families of $6n$-dimensional representations obtained from the quotient singularity $\C^2/\Z_6$.
\end{abstract}

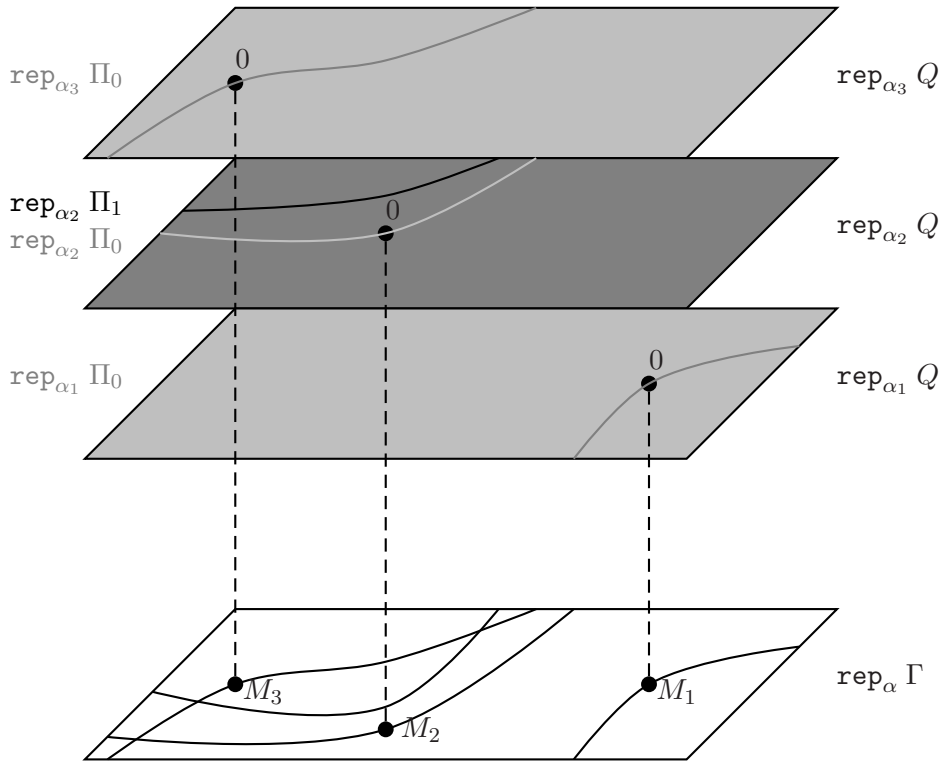
\[
\begin{pspicture}(-1,0)(15,10)
\pspolygon[fillstyle=solid,fillcolor=lightgray](0,8)(2,10)(10,10)(8,8)
\pspolygon[fillstyle=solid,fillcolor=gray](0,6)(2,8)(10,8)(8,6)
\pspolygon[fillstyle=solid,fillcolor=lightgray](0,4)(2,6)(10,6)(8,4)
\pspolygon[fillstyle=solid,fillcolor=white](0,0)(2,2)(10,2)(8,0)
\cnode*(2,9){3pt}{A}
\cnode*(4,7){3pt}{B}
\cnode*(7.5,5){3pt}{C}
\cnode*(2,1){3pt}{A1}
\cnode*(4,.4){3pt}{B1}
\cnode*(7.5,1){3pt}{C1}
\ncline[linestyle=dashed]{A}{A1}
\ncline[linestyle=dashed]{B}{B1}
\ncline[linestyle=dashed]{C}{C1}
\pscurve[linecolor=gray](0.3,8)(2,9)(4,9.3)(6,10)
\pscurve(.3,0)(2,1)(4,1.3)(6,2)
\pscurve[linecolor=lightgray](1,7)(4,7)(6,8)
\pscurve(1.3,7.3)(4,7.5)(5.5,8)
\pscurve(.9,.9)(4,.7)(5.5,2)
\pscurve(.3,.3)(4,.4)(6.5,2)
\pscurve[linecolor=gray](6.5,4)(7.5,5)(9.5,5.5)
\pscurve(6.5,0)(7.5,1)(9.5,1.5)
\put(10,1){$\wis{rep}_{\alpha}~\Gamma$}
\put(10,5){$\wis{rep}_{\alpha_1}~Q$}
\put(10,7){$\wis{rep}_{\alpha_2}~Q$}
\put(10,9){$\wis{rep}_{\alpha_3}~Q$}
\put(-1,9){\gray $\wis{rep}_{\alpha_3}~\Pi_0$}
\put(-1,6.8){\gray $\wis{rep}_{\alpha_2}~\Pi_0$}
\put(-1,7.3){\black $\wis{rep}_{\alpha_2}~\Pi_1$}
\put(-1,5){\gray $\wis{rep}_{\alpha_1}~\Pi_0$}
\put(2,9.2){$0$}
\put(4,7.2){$0$}
\put(7.5,5.2){$0$}
\put(7.6,.8){$M_1$}
\put(2.1,.8){$M_3$}
\put(4.2,.3){$M_2$}
\end{pspicture}
\]

\section{Introduction}

The classification of semi-simple representations of the modular group $\Gamma = SL_2(\Z)$ is of crucial importance to many branches of mathematics. In particular, to topological quantum field theory via modular tensor categories \cite{BakalovKirillov} and to knot theory via its relation with the third braid group $B_3$, see for example \cite{Westbury2}. Still, very little seems to be known on this problem except for representations of dimensions $\leq 5$ by \cite{TubaWenzl}.

In this paper we will use techniques from non-commutative geometry to get a grip on (nice families of) $n$-dimensional representations $\wis{rep}_n~\Gamma$ of $\Gamma$ for arbitrary $n$. The starting point is that the modular group is an amalgamated free product 
\[
SL_2(\Z) = \Gamma = \Z_4 \ast_{\Z_2} \Z_6 \]
(see for example \cite{Serre}) and hence that its group algebra $\C \Gamma$ is a formally smooth algebra in the sense of \cite{CuntzQuillen} or \cite{KontRos}. As such we can use the general techniques outlined in \cite{LBnagatn} to study the finite dimensional representations of $\C \Gamma$ via specific quiver-representations encoding the \'etale local structure of $\wis{rep}_n~\Gamma$ in the neighborhood of a semi-simple representation.

An approach to the representation theory of the {\em projective} modular group $\overline{\Gamma} = PSL_2(\Z)$ via quiver-representations was already found by Bruce Westbury in \cite{Westbury1} and can be modified to be applicable for $\Gamma$. The idea is to study the restrictions of an $n$-dimensional representation $M$ to the cyclic subgroups $\Z_2,\Z_4$ and $\Z_6$ of $\Gamma$ and decompose these in eigenspaces. In this way one associates to $M$ a representation of the bipartite quiver $\Lambda$
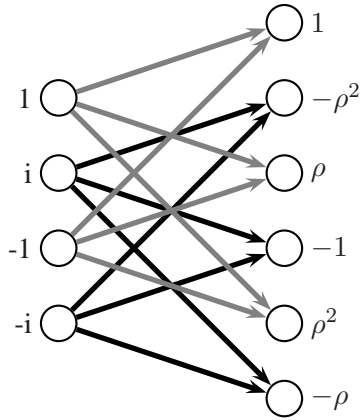
\[
\begin{pspicture}(-.5,-.5)(3.5,5.5)
\cnode(0,1){.25}{A}
\cnode(0,2){.25}{B}
\cnode(0,3){.25}{C}
\cnode(0,4){.25}{D}
\cnode(3,0){.25}{E}
\cnode(3,1){.25}{F}
\cnode(3,2){.25}{G}
\cnode(3,3){.25}{H}
\cnode(3,4){.25}{I}
\cnode(3,5){.25}{J}
\uput{10pt}[l](0,1){-i}
\uput{10pt}[l](0,2){-1}
\uput{10pt}[l](0,3){i}
\uput{10pt}[l](0,4){1}
\uput{10pt}[r](3,0){$-\rho$}
\uput{10pt}[r](3,1){$\rho^2$}
\uput{10pt}[r](3,2){$-1$}
\uput{10pt}[r](3,3){$\rho$}
\uput{10pt}[r](3,4){$-\rho^2$}
\uput{10pt}[r](3,5){$1$}
\ncline[linecolor=black,linewidth=2pt]{->}{A}{E}
\ncline[linecolor=black,linewidth=2pt]{->}{A}{G}
\ncline[linecolor=black,linewidth=2pt]{->}{A}{I}
\ncline[linecolor=black,linewidth=2pt]{->}{C}{I}
\ncline[linecolor=black,linewidth=2pt]{->}{C}{G}
\ncline[linecolor=black,linewidth=2pt]{->}{C}{E}
\ncline[linecolor=gray,linewidth=2pt]{->}{B}{F}
\ncline[linecolor=gray,linewidth=2pt]{->}{B}{H}
\ncline[linecolor=gray,linewidth=2pt]{->}{B}{J}
\ncline[linecolor=gray,linewidth=2pt]{->}{D}{F}
\ncline[linecolor=gray,linewidth=2pt]{->}{D}{H}
\ncline[linecolor=gray,linewidth=2pt]{->}{D}{J}
\end{pspicture}
\]
Here, the left vertices correspond to the eigenvalues of $\Z_4$, those on the right to those of $\Z_6$ ($\rho$ is a primitive third root of unity) and the arrows describe the base-change matrix
\[
\psi~:~M \downarrow_{\Z_4} \rTo M \downarrow_{\Z_6} \]
which is $\Z_2$-invariant. The corresponding dimension vector $\alpha = (a_1,\hdots,a_4 ; b_1,\hdots,b_6)$ describes the characters of these restrictions and satisfies $\sum_i a_i = n = \sum_j b_j$. Moreover, the quiver-representation $R_M$ of $\wis{rep}_{\alpha}~\Lambda$ corresponding to $M$ has the property that the map $\psi$ is invertible whence $R_M$ must be $\theta$-stable for the stability structure $\theta = (-1,\hdots,-1;1,\hdots,1)$. In this way one can decompose $\wis{rep}_n~\Gamma$ as a disjoint union of connected components $\wis{rep}_{\alpha}~\Gamma$ each of which can be identified with an affine open subset of $\wis{moduli}^{\theta}_{\alpha}~\Lambda$, the moduli space of $\theta$-semistable representations of dimension vector $\alpha$ of $\Lambda$, see for example \cite{King}.

It turns out that moduli spaces of quiver-representations are much harder to study than algebraic quotient varieties of quiver-representations. In our previous paper \cite{AdriLB} we used local quivers to reduce the \'etale local structure of these moduli spaces to ordinary quotient varieties and this technique could be applied here, too. However, we choose a different route inspired by \cite{LBOneQuiver} where a quiver-setting $(Q,\mu)$ is associated to any formally smooth algebra $A$, the conjectural explanation of which is that there should be a non-commutative \'etale isomorphism between $A$ and an algebra $B$ which is Morita equivalent to the path algebra $\C Q$ with the Morita equivalence determined by the dimension vector $\mu$.
If we apply this procedure to $\C \Gamma$ we obtain the quiver $Q$
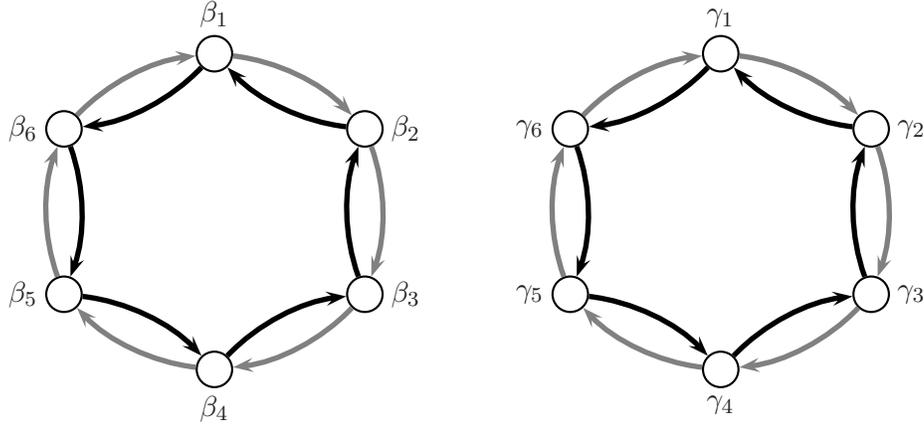
\[
\begin{pspicture}(-.5,-.5)(4.5,4.7)
\cnode(2,4.2){.25}{A}
\cnode(4,3.2){.25}{B}
\cnode(4,1){.25}{C}
\cnode(2,0){.25}{D}
\cnode(0,1){.25}{E}
\cnode(0,3.2){.25}{F}
\uput{10pt}[u](2,4.2){$\beta_1$}
\uput{10pt}[r](4,3.2){$\beta_2$}
\uput{10pt}[r](4,1){$\beta_3$}
\uput{10pt}[d](2,0){$\beta_4$}
\uput{10pt}[l](0,1){$\beta_5$}
\uput{10pt}[l](0,3.2){$\beta_6$}
\ncarc[arcangle=20,linecolor=gray,linewidth=2pt]{->}{A}{B}
\ncarc[arcangle=20,linecolor=gray,linewidth=2pt]{->}{B}{C}
\ncarc[arcangle=20,linecolor=gray,linewidth=2pt]{->}{C}{D}
\ncarc[arcangle=20,linecolor=gray,linewidth=2pt]{->}{D}{E}
\ncarc[arcangle=20,linecolor=gray,linewidth=2pt]{->}{E}{F}
\ncarc[arcangle=20,linecolor=gray,linewidth=2pt]{->}{F}{A}
\ncarc[arcangle=20,linecolor=black,linewidth=2pt]{->}{A}{F}
\ncarc[arcangle=20,linecolor=black,linewidth=2pt]{->}{F}{E}
\ncarc[arcangle=20,linecolor=black,linewidth=2pt]{->}{E}{D}
\ncarc[arcangle=20,linecolor=black,linewidth=2pt]{->}{D}{C}
\ncarc[arcangle=20,linecolor=black,linewidth=2pt]{->}{C}{B}
\ncarc[arcangle=20,linecolor=black,linewidth=2pt]{->}{B}{A}
\end{pspicture}~\qquad~\qquad
\begin{pspicture}(-.5,-.5)(4.5,4.7)
\cnode(2,4.2){.25}{A}
\cnode(4,3.2){.25}{B}
\cnode(4,1){.25}{C}
\cnode(2,0){.25}{D}
\cnode(0,1){.25}{E}
\cnode(0,3.2){.25}{F}
\uput{10pt}[u](2,4.2){$\gamma_1$}
\uput{10pt}[r](4,3.2){$\gamma_2$}
\uput{10pt}[r](4,1){$\gamma_3$}
\uput{10pt}[d](2,0){$\gamma_4$}
\uput{10pt}[l](0,1){$\gamma_5$}
\uput{10pt}[l](0,3.2){$\gamma_6$}
\ncarc[arcangle=20,linecolor=gray,linewidth=2pt]{->}{A}{B}
\ncarc[arcangle=20,linecolor=gray,linewidth=2pt]{->}{B}{C}
\ncarc[arcangle=20,linecolor=gray,linewidth=2pt]{->}{C}{D}
\ncarc[arcangle=20,linecolor=gray,linewidth=2pt]{->}{D}{E}
\ncarc[arcangle=20,linecolor=gray,linewidth=2pt]{->}{E}{F}
\ncarc[arcangle=20,linecolor=gray,linewidth=2pt]{->}{F}{A}
\ncarc[arcangle=20,linecolor=black,linewidth=2pt]{->}{A}{F}
\ncarc[arcangle=20,linecolor=black,linewidth=2pt]{->}{F}{E}
\ncarc[arcangle=20,linecolor=black,linewidth=2pt]{->}{E}{D}
\ncarc[arcangle=20,linecolor=black,linewidth=2pt]{->}{D}{C}
\ncarc[arcangle=20,linecolor=black,linewidth=2pt]{->}{C}{B}
\ncarc[arcangle=20,linecolor=black,linewidth=2pt]{->}{B}{A}
\end{pspicture}
\]
where the vertices correspond to the twelve one-dimensional simple representations $\{ S_1,\hdots,S_{12} \}$ of $\Gamma$. We obtain an algebra map
\[
\phi~:~\C \Gamma \rTo \C Q \]
which should be non-commutative \'etale in a yet to be developed non-commutative Grothendieck topology (we make some comments on this in the final section). We have an induced map
\[
\phi^*~:~\wis{rep}~Q \rTo \wis{rep}~\Gamma \]
the properties of which we will now describe (see the frontispiece as illustration). Each component $\wis{rep}_{\alpha}~\Gamma$ of $n$-dimensional representations of $\Gamma$ contains several semi-simple representations which are made up of the twelve simples
\[
M_i = S_1^{\oplus e_1(i)} \oplus \hdots \oplus S_{12}^{\oplus e_{12}(i)} \]
which determine a dimension vector $\alpha_i = (e_1(i),\hdots,e_{12}(i))$ of total dimension $n$ for $Q$. By the results of \cite{AdriLB} we know that the restriction of $\phi^*$ to $\wis{rep}_{\alpha_i}~Q$ determines a $GL_n$-equivariant \'etale map between a neighborhood of the orbit of $M_i$ in $\wis{rep}_{\alpha}~\Gamma$ and a neighborhood of the zero-representation in $\wis{rep}_{\alpha_i}~Q$. As such we also have an \'etale isomorphism between $\wis{iss}_{\alpha}~\Gamma$ (the variety parametrizing isoclasses of semi-simple $\Gamma$-representations in $\wis{rep}_{\alpha}~\Gamma$) near the point $[M_i]$ and the quiver quotient-variety $\wis{iss}_{\alpha_i}~Q$ near the point $[0]$.

The quiver $Q$ carries a symplectic structure, whence we can apply the results of \cite{LBBocklandt} to define symplectic flows (similar to the Calogero-Moser flows studied by G. Wilson in \cite{Wilson}) on $\wis{rep}_{\alpha_i}~Q$ and the quotient varieties $\wis{iss}_{\alpha_i}~Q$ which allow us to transport local information near $[0]$ to an arbitrary semi-simple $\alpha_i$-dimensional representation of $Q$. In this way we can describe (at least in principle) most of the semi-simple $n$-dimensional representations of $\Gamma$. We do not know whether all $\wis{iss}_{\alpha_i}~Q$ for a fixed $\alpha$ cover the whole of $\wis{iss}_{\alpha}~\Gamma$ but conjecture that they do.

To construct semi-simple $\Gamma$-representations with specific properties (see for example \cite{Westbury2} for the construction of knot invariants) we need to have rational families in $\wis{iss}_{\alpha}~\Gamma$. Here too, the symplectic structure of $\wis{rep}_{\alpha_i}~Q$ comes in handy. The knowledgeable reader will have noticed that each components of $Q$ is the quiver which Kronheimer \cite{Kronheimer} and others use to deform the quotient singularity $\C^2/\Z_6$, see for example \cite{Slodowy} and \cite{Crawley} for more details. This allows us (via the map $\phi^*$) to give very explicit $2n$-parameter families of simple $\Gamma$-representations of dimension $6n$. We can even use this to construct several nice families of semi-simple representations in $\wis{rep}_{\alpha_i}~Q$ (the 'curves' in the frontispiece) which give us interesting filaments of semi-simple $\Gamma$-representations. We will return to the detailed study of these families elsewhere.

\par \vskip 3mm

Much of the strategy outlined above is applicable to more general situations and we formulate our results accordingly. In section 2 we outline the procedure to associate the \'etale type quiver-setting $(Q,\mu)$ to any amalgamated coproduct $S_1 \ast_{S} S_2$ of finite dimensional semi-simple $\C$-algebras $S,S_1$ and $S_2$. In section 3 we give several examples, including the modular group example as well as coproducts of matrix-algebras. In section 4 we will recall results from \cite{LBBocklandt} on non-commutative symplectic geometry and compute some symplectic flows. We also observe that the fact that the \'etale type quiver is symplectic is valid in great generality. In section 5 we make the algebra map $\phi$ explicit and outline the method to find $\Gamma$-representations from the quotient singularity $\C^2/\Z_6$. In the final section we make some comments on non-commutative \'etale algebra maps.

\section{Amalgamating semi-simple algebras}

A {\em formally smooth algebra} $A$ has the lifting property with respect to nilpotent ideals (generalizing Grothendieck's characterization of commutative regular algebras), that is, any diagram of associative $\C$-algebras with unit
\[
\begin{diagram}
T & \rOnto^{\pi} & T/I \\
& \luDotsto_{\tilde{\phi}} & \uTo^{\phi} \\
& & A
\end{diagram}
\]
can be completed with an algebra morphism $\phi$ whenever $I$ is a nilpotent ideal of $T$. These algebras have been studied by J. Cuntz and D. Quillen in \cite{CuntzQuillen} where they are called
{\em quasi-free algebras}. Affine formally smooth algebras are thought to be the coordinate rings of {\em non-commutative manifolds}, see for example \cite{KontRos}.  Much as a manifold is locally diffeomorphic to an affine space, a formally smooth algebra $A$ is conjectured to have a non-commutative \'etale extension $B$ which is the coordinate ring of a {\em non-commutative affine space} which are (conjecturally) rings Morita equivalent to path algebras of quivers. In \cite{LBOneQuiver} a procedure was given to associate to a formally smooth affine algebra $A$ a {\em quiver setting} $(Q_A,\alpha_A)$ where $Q_A$ is a finite quiver and $\alpha_A$ a dimension vector determining the Morita equivalent ring $B$ associated to $A$. We quickly run through the construction and refer to \cite{LBOneQuiver} for more details.

Because $A$ is formally smooth, the scheme of $n$-dimensional representations $\wis{rep}_n~A$ is a smooth reduced variety, possibly having several connected components
\[
\wis{rep}_n~A = \bigsqcup_{| \alpha |=n} \wis{rep}_{\alpha}~A \]
which we give formal symbols $\alpha$ and we say that $\alpha$ is a dimension vector of dimension $n$ and denote $| \alpha | = n$. The direct sum of representations induces a semigroup structure on the collection of all $\alpha$ for all $n \in \N$. The vertices of $Q_A$ are in one-to-one correspondence with the semigroup generators of this semigroup (they are uniquely determined as they correspond to those components which are principal $PGL_n$-fibrations) say $\{ \beta_1,\hdots,\beta_k \}$. The number of directed arrows in $Q_A$ from vertex $i$ to vertex $j$ is given by $\wis{ext}(\beta_i,\beta_j)$ which is the minimal dimension of $Ext_A^1(S,S')$ with $S \in \wis{rep}_{\beta_i}~A$ and $S' \in \wis{rep}_{\beta_j}~A$. The dimension vector $\alpha_A$ has as its $i$-th component $| \beta_i |$. Conjecturally, the non-commutative \'etale type of $A$ is determined by the algebra which is Morita-equivalent (determined by $\alpha_A$) with the path algebra $\C Q_A$. We will explain this construction in a large class of examples : {\em amalgamated semi-simple algebras}.

Consider three finite-dimensional semi-simple $\C$-algebras $S$, $S_1$ and $S_2$ with decompositions
\[
S = M_{k_1}(\C) \oplus \hdots \oplus M_{k_r}(\C) \]
\[
S_1 = M_{l_1}(\C) \oplus \hdots \oplus M_{l_s}(\C) \]
\[
S_2 = M_{m_1}(\C) \oplus \hdots \oplus M_{m_t}(\C) \]
and assume we have embeddings $i_1~:~S \rInto S_1$ and $i_2~:~S \rInto S_2$ which are fully determined by their {\em Bratelli diagrams}. Recall that these are labeled full bipartite graphs on $(r,s)$ (resp. $(r,t)$) vertices
where the edge connecting the $i$-th vertex on the left to the $j$-th vertex on the right is labeled with a natural number $a_{ij} \geq 0$ (resp. $b_{ij} \geq 0$) satisfying the following numerical restrictions
\[
l_j = \sum_{i=1}^r a_{ij} k_i \qquad \text{and} \qquad m_i = \sum_{j=1}^r b_{ij}k_i \]
These number encode the {\em restriction data} : let $U_1,\hdots,U_r$ be the simples of $S$ (of $\C$-dimensions $k_1,\hdots,k_r$) and similarly $V_1,\hdots,V_s$ (resp. $W_1,\hdots,W_t$) the simples of $S_1$ (resp. of $S_2$) which have $\C$-dimensions $l_1,\hdots,l_r$ (resp. $m_1,\hdots,m_s$). As an $S$-module, the $V_i$ and $W_i$ decompose in a direct sum of the $U_j$ and the numbers $a_{ij}$ resp. $b_{ij}$ denote the multiplicities, that is
\[
V_j \downarrow_{S} \simeq U_1^{\oplus a_{1j}} \oplus \hdots \oplus U_r^{\oplus a_{rj}} \quad \text{and} \quad
W_j \downarrow_{S} \simeq U_1^{\oplus b_{1j}} \oplus \hdots \oplus U_r^{\oplus b_{rj}}
\]
From these decompositions and Schur's lemma it follows that
\[
Hom_{S}(V_i,W_j) = \sum_{k=1}^r a_{ki}b_{kj} = n_{ij} \]
We are interested in the {\em amalgamated coproduct algebra}
\[
A = S_1 \ast_S S_2 \]
which represents the functor $Hom_{S-alg}(S_1,-) \times Hom_{S-alg}(S_2,-)$ in the category of all  $S$-algebras. In \cite{Schofieldbook} it was proved that $A$ is hereditary. In fact, 

\begin{proposition} For $S,S_1,S_2$ finite-dimensional semi-simple $\C$-algebras as above, the amalgamated coproduct
\[
A = S_1 \ast_S S_2 \]
is a formally smooth algebra. In fact,  if $S_1$ and $S_2$ are formally smooth and $S$ is finite-dimensional semi-simple, the same result holds. 
\end{proposition}

\begin{proof} Let $T$ be a $\C$-algebra and $I \triangleleft T$ a nilpotent ideal and assume we have a morphism $\phi$
\[
\begin{diagram}
S_1 & \rTo^{\alpha} & S_1 \ast_S S_2 & \rTo^{\phi} & T/I \\
\uInto & & \uTo^{\beta} & & \uOnto^{\pi} \\
S & \rInto & S_2 & & T 
\end{diagram}
\]
Because $S_1$ (resp. $S_2$) is formally smooth there is a lift $\tilde{\alpha}~:~S_1 \rTo T$ (resp. $\tilde{\beta}~:~S_2 \rTo T$) compatible with $\phi \circ \alpha$ (resp. with $\phi \circ \beta$). By
\cite[Prop. 6.1]{CuntzQuillen} there is a unit $t \in T^*$ such that $\pi(t)=1$ satisfying
\[
\tilde{\alpha} \mid S = t^{-1} (\tilde{\beta} \mid S) t \]
But then, $\tilde{\alpha} \ast t^{-1} \tilde{\beta} t$ is a required lift from $S_1 \ast_S S_2$ to $T$.
\end{proof}

We want to determine the connected components of $\wis{rep}_n~A$. Let $M$ be an $n$-dimensional representation of $A$, then restricting $M$ to $S_1$ and $S_2$ gives us decompositions
\[
M \downarrow_{S_1} \simeq V_1^{\oplus p_1} \oplus \hdots \oplus V_s^{\oplus p_s} \quad \text{and} \quad
M \downarrow_{S_2} \simeq W_1^{\oplus q_1} \oplus \hdots \oplus W_t^{\oplus q_t} \]
Clearly, the $s+t$-tuple of natural numbers $\alpha = (p_1,\hdots,p_s;q_1,\hdots;q_t)$ satisfies the numerical restriction
\[
n = l_1p_1+\hdots+l_sp_s = m_1q_1+\hdots+m_tq_t \]
Choose a $\C$-basis in $\C^n = M \downarrow_{S_1}$ compatible with the decomposition, say
$\mathcal{B} = \{ b_1,\hdots,b_n \}$ and a $\C$-basis $\mathcal{B'} = \{ b'_1,\hdots,b'_n \}$ of
$\C^n = M \downarrow_{S_2}$ compatible with the second decomposition. Because $M$ is a
$S_1 \ast_S S_2$-representation, the base-change map $\mathcal{B} \rTo^{\psi} \mathcal{B'}$ must be an invertible element of
\[
Hom_S(M \downarrow_{S_1},M \downarrow_{S_2}) = \oplus_{i=1}^s \oplus_{j=1}^t M_{q_j \times p_i}(Hom_S(V_i,W_j)) \]
and hence $\psi$ determines a representation of dimension vector $\alpha$ of the bipartite quiver $\Lambda$ on $(s,t)$ vertices (see figure~1)
\begin{figure} \label{quiver}
\[
\begin{pspicture}(0,-2)(3,7)
\cnode(0,4.5){.25}{C}
\cnode(0,2.8){.25}{B}
\cnode(0,.5){.25}{A}
\put(0,1.5){$\vdots$}
\cnode(4,7){.25}{D}
\cnode(4,5){.25}{E}
\cnode(4,-2){.25}{F}
\put(4,1.5){$\vdots$}
\ncline[linecolor=gray,linewidth=2pt]{->}{C}{D} \ncput*{\tiny{$n_{11}$}}
\ncline[linecolor=gray,linewidth=2pt]{->}{C}{E}
\ncline[linecolor=gray,linewidth=2pt]{->}{C}{F} 
\ncline[linecolor=darkgray,linewidth=2pt]{->}{B}{D} 
\ncline[linecolor=darkgray,linewidth=2pt]{->}{B}{E} 
\ncline[linecolor=darkgray,linewidth=2pt]{->}{B}{F}
\ncline[linecolor=black,linewidth=2pt]{->}{A}{D} 
\ncline[linecolor=black,linewidth=2pt]{->}{A}{E}
\ncline[linecolor=black,linewidth=2pt]{->}{A}{F} \ncput*{\tiny{$n_{st}$}}
\end{pspicture}
\]
\caption{The bipartite quiver $\Lambda$.}
\end{figure}
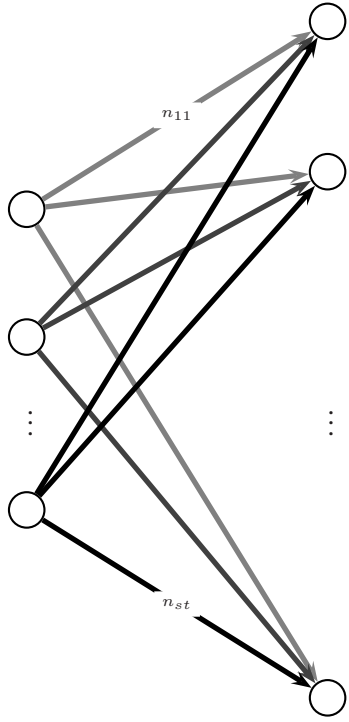
having $n_{ij} = \sum_{k=1}^r a_{ki}b_{kj}$ arrows pointing from the $i$-th vertex on the left to the $j$-th vertex on the right. Clearly, these numbers of arrows must be so that there is a representation in $V_M \in \wis{rep}_{\alpha}~\Lambda$ such that the $n \times n$ matrix describing $\psi$ (and hence $M$) is invertible. This implies that the $\Lambda$-representation $V_M$ is {\em $\theta$-semi-stable} where $\theta=(-l_1,\hdots,-l_s;m_1,\hdots,m_t)$ and $\theta$-semi-stability means that there is no proper subrepresentation $W$ of $V_M$ with dimension vector $\beta$ such that $\theta.\beta < 0$, see for example \cite{AdriLB} for more details. Summarizing this argument, we have

\begin{proposition} Let $M \in \wis{rep}_n~S_1 \ast_S S_2$, then there is a point $M'$ in the $GL_n$-orbit $\Oscr(M)$ of $M$ such that $M'$ is determined by a representation
\[
V_{M'} \in \wis{rep}_{\alpha}~\Lambda \]
with $\alpha=(p_1,\hdots.p_s;q_1,\hdots,q_t) \in \N^{s+t}$ such that
\[
n = \sum_{i=1}^s l_ip_i = \sum_{j=1}^t m_jq_j \]
If $M"$ is another such point in the orbit $\Oscr(M)$, then $V_{M'}$ and $V_{M"}$ are isomorphic as $\Lambda$-representations, that is, there is a basechange in $GL(\alpha)$ mapping $V_{M'}$ to $V_{M"}$.
\end{proposition}

As a consequence we see that $\wis{rep}_n~S_1 \ast_S S_2$ has as many connected components as there are dimension vectors $\alpha = (p_1,\hdots,p_s;q_1,\hdots,q_t)$ such that (1) $\sum l_ip_i = n = \sum_j m_jq_j$ and (2)  $\wis{rep}_{\alpha}~Q$ contains $\theta$-semi-stable representations. We will denote the set of all such dimension vectors with $\wis{comp} \subset \N^{s+t}$. $\wis{comp}$ is a semigroup having a finite set of semigroup-generators $\{ \beta_1,\hdots,\beta_k \}$ determined by the property that $\wis{rep}_{\beta_i}~\Lambda$ consists entirely of $\theta$-stable representations (those for which every proper subrepresentation has a dimension vector $\beta$ satisfying $\theta.\beta > 0$).

The {\em Euler-form} $\chi_{\Lambda}$ is the bilinear form on $\Z^{s+t}$ determined by the integral matrix 
\[
\chi_{\Lambda} = \begin{bmatrix}
\begin{array}{ccc|ccc}
1 & & & -n_{11} & \hdots & -n_{1t} \\
& \ddots & & \vdots & & \vdots \\
& & 1 & -n_{s1} & \hdots & -n_{st} \\
\hline 
0 & \hdots & 0 & 1 & & \\
\vdots & & \vdots & & \ddots & \\
0 & \hdots & 0 & & & 1
\end{array}
\end{bmatrix}
\]
Recall that if $V$ resp. $W$ are representations of $\Lambda$ of dimension vector $\beta$ resp. $\gamma$, then we have the numerical relation
\[
dim_{\C}~Hom_{\Lambda}(V,W) - dim_{\C}~Ext^1_{\Lambda}(V,W) = \chi_{\Lambda}(\beta,\gamma)
\]
which allows us to compute the values of $\wis{ext}(\beta_i,\beta_j)$ observing that for $\theta$-stable representations $V$ and $W$ we have that $dim_{\C}~Hom_{\Lambda}(V,W) = \delta_{VW}$. We have now all the required material to determine the quiver setting associated to the amalgamated coproduct.

\begin{theorem} With notation as introduced before, the quiver setting $(Q_A,\alpha_A)$ associated to the amalgamated coproduct of semi-simple algebras
\[
A = S_1 \ast_S S_2 \]
is given by the quiver $Q_A$ on $k$ vertices where $\{ \beta_1,\hdots,\beta_k \}$ is the set of semigroup-generators of the semigroup $\wis{comp} \subset \N^{s+t}$ consisting of all $\theta$-semi-stable dimension vectors of $\Lambda$ and where the number of directed arrows from vertex $i$ to vertex $j$ is given by the number
\[
\delta_{ij} - \chi_{\Lambda}(\beta_i,\beta_j) \]
The dimension vector $\alpha_A = (n_1,\hdots,n_k) \in \N^k$ is determined by 
\[
n_i = \sum_{j=1}^s l_j p_j(i) \]
if $\beta_i = (p_1(i),\hdots,p_s(i);q_1(i),\hdots,q_t(i))$.
\end{theorem}

\section{The \'etale type of the modular group.}

In this section we will give several examples of the foregoing construction with emphasis on the important case of the group algebra $\C SL_2(\Z)$ of the modular group. To warm-up let us begin with a tame example.

\begin{example}[The infinite dihedral group] As $D_{\infty} = \langle s,t~|~s^2=1=t^2 \rangle \simeq \Z_2 \ast \Z_2$ we have that
\[
\C D_{\infty} \simeq (\C \times \C) \ast (\C \times \C) \]
with the factors on both sides corresponding to the eigenvalues $\pm 1$. The Bratelli diagrams for the embedding of $\C$ in these factors are
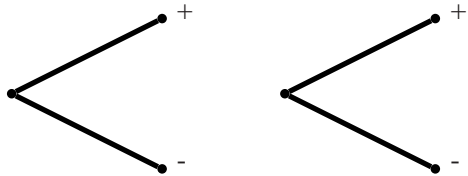
\[
\begin{pspicture}(0,-1)(2,1)
\dotnode(0,0){A}
\dotnode(2,1){B}
\put(2.2,1){+}
\dotnode(2,-1){C}
\put(2.2,-1){-}
\ncline[linecolor=black,linewidth=2pt]{A}{B} 
\ncline[linecolor=black,linewidth=2pt]{A}{C}
\end{pspicture} \qquad~\qquad
\begin{pspicture}(0,-1)(2,1)
\dotnode(0,0){A}
\dotnode(2,1){B}
\put(2.2,1){+}
\dotnode(2,-1){C}
\put(2.2,-1){-}
\ncline[linecolor=black,linewidth=2pt]{A}{B} 
\ncline[linecolor=black,linewidth=2pt]{A}{C}
\end{pspicture}
\]
whence the quiver $\Lambda$ is given by
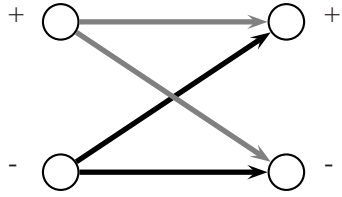
\[
\begin{pspicture}(-1,-.5)(4,2.5)
\cnode(0,0){.25}{A}
\cnode(0,2){.25}{B}
\cnode(3,0){.25}{C}
\cnode(3,2){.25}{D}
\put(-.7,0){-}
\put(-.7,2){+}
\put(3.5,0){-}
\put(3.5,2){+}
\ncline[linecolor=black,linewidth=2pt]{->}{A}{C}
\ncline[linecolor=black,linewidth=2pt]{->}{A}{D}
\ncline[linecolor=gray,linewidth=2pt]{->}{B}{C}
\ncline[linecolor=gray,linewidth=2pt]{->}{B}{D}
\end{pspicture}
\]
The character $\theta=(-1,-1;1,1)$ whence $\wis{comp}$ is the semigroup of $\N^4$ consisting of dimension vectors $(a,b;c,d)$ such that $a+b=c+d$. This semigroup has $4$ generators
\[
\begin{cases}
++=\beta_1=(1,0;1,0) \\ +-=\beta_2=(1,0;0,1) \\ -+=\beta_3=(0,1;1,0) \\ --=\beta_4=(0,1;0,1)
\end{cases}
\]
The quiver setting $(Q,\alpha)$ is given by $\alpha=(1,1,1,1)$ with $Q$ the quiver
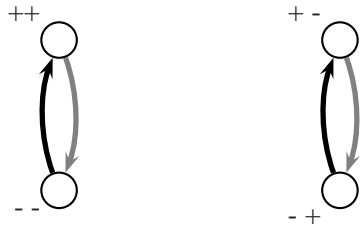
\[
\begin{pspicture}(0,-.5)(2,2.5)
\cnode(1,0){.25}{A}
\uput{10pt}[dl](1,0){- -}
\cnode(1,2){.25}{B}
\uput{10pt}[ul](1,2){++}
\ncarc[arcangle=20,linecolor=gray,linewidth=2pt]{->}{B}{A}
\ncarc[arcangle=20,linecolor=black,linewidth=2pt]{->}{A}{B}
\end{pspicture}~\qquad~\qquad
\begin{pspicture}(0,-.5)(2,2.5)
\cnode(1,0){.25}{A}
\uput{10pt}[dl](1,0){- +}
\cnode(1,2){.25}{B}
\uput{10pt}[ul](1,2){+ -}
\ncarc[arcangle=20,linecolor=gray,linewidth=2pt]{->}{B}{A}
\ncarc[arcangle=20,linecolor=black,linewidth=2pt]{->}{A}{B}
\end{pspicture}
\]
\end{example}

The study of isomorphism classes of semi-simple representations of the {\em modular group}
$\Gamma = SL_2(\Z)$ is of fundamental importance in the study of modular tensor categories and hence in topological quantum field theory, see for example \cite{BakalovKirillov}. Still, surprisingly little is known about such representations. To the best of our knowledge the best results are contained in \cite{TubaWenzl} and classify the simple representations in dimension $\leq 5$.

In order to fix notation let us recall the relation between $\Gamma$ and the {\em three string braid group} which is generated by the elementary braids
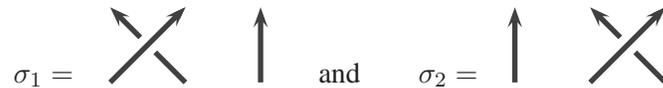
\[
\sigma_1 =\quad
\begin{pspicture}(0,0)(2,1)
\psline[linecolor=darkgray,linewidth=2pt]{->}(1,0)(0,1)
\psline[linecolor=darkgray,linewidth=2pt]{->}(2,0)(2,1)
\psline[linecolor=darkgray,linewidth=2pt,border=3pt]{->}(0,0)(1,1)
\end{pspicture} \qquad \text{and} \qquad
\sigma_2 =\quad
\begin{pspicture}(0,0)(2,1)
\psline[linecolor=darkgray,linewidth=2pt]{->}(2,0)(1,1)
\psline[linecolor=darkgray,linewidth=2pt]{->}(0,0)(0,1)
\psline[linecolor=darkgray,linewidth=2pt,border=3pt]{->}(1,0)(2,1)
\end{pspicture}
\]
which satisfy the braid group relation : $\sigma_1 \sigma_2 \sigma_1 = \sigma_2 \sigma_1 \sigma_2$
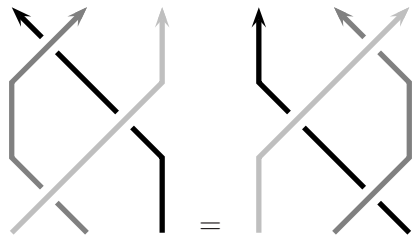
\[
\begin{pspicture}(0,0)(2,3)
\psline[linecolor=black,linewidth=2pt]{->}(2,0)(2,1)(0,3)
\psline[border=3pt,linecolor=gray,linewidth=2pt]{->}(1,0)(0,1)(0,2)(1,3)
\psline[border=3pt,linecolor=lightgray,linewidth=2pt]{->}(0,0)(1,1)(2,2)(2,3)
\end{pspicture} \quad = \quad
\begin{pspicture}(0,0)(2,3)
\psline[linecolor=black,linewidth=2pt]{->}(2,0)(0,2)(0,3)
\psline[border=3pt,linecolor=gray,linewidth=2pt]{->}(1,0)(2,1)(2,2)(1,3)
\psline[border=3pt,linecolor=lightgray,linewidth=2pt]{->}(0,0)(0,1)(2,3)
\end{pspicture}
\]
If we define $s = \sigma_1 \sigma_2 \sigma_1$ and $t= \sigma_1 \sigma_2$, $B_3$ has another presentation
\[
 B_3 = \langle s,t~:~s^2=t^3 \rangle \]
 From this description it is clear that $c = s^2 = t^3$ is a central element and that the corresponding quotient group
\[
\overline{B}_3 = B_3/\langle c \rangle \simeq \Z_2 \ast \Z_3 \]
is the free product of cyclic groups of order $2$ and $3$. It is well-known that the braid group $B_3$ maps onto  $SL_2(\Z)$ via the map
\[
\sigma_1 \mapsto \begin{bmatrix} 1 & 1 \\ 0 & 1 \end{bmatrix} \qquad \text{and} \qquad
\sigma_2 \mapsto \begin{bmatrix} 1 & 0 \\ -1 & 1 \end{bmatrix}
\]
which maps the generators $s$ and $t$, respectively the central element $c$, to
\[
s \mapsto \begin{bmatrix} 0 & 1 \\ -1 & 0 \end{bmatrix} \qquad \text{and} \qquad
t \mapsto \begin{bmatrix} 1 & 1 \\ -1 & 0 \end{bmatrix} \quad \text{resp.} \qquad c \mapsto \begin{bmatrix} -1 & 0 \\ 0 & -1 \end{bmatrix}  \]
Therefore, the {\em projective modular group} $\overline{\Gamma} = PSL_2(\Z) \simeq \overline{B}_3 \simeq \Z_2 \ast \Z_3$ whereas the full modular group has the structure of an {\em amalgamated free product}
\[
SL_2(\Z) = \Gamma \simeq \Z_4 \ast_{\Z_2} \Z_6 \]
see for example \cite[4.2(c)]{Serre}.  In this paper we will initiate the study of the representation theory of these groups from a non-commutative geometry perspective.

\begin{example}[The projective modular group] As $\overline{\Gamma} \simeq \Z_2 \ast \Z_3$ we have that
\[
\C~\Gamma \simeq (\C \times \C) \ast (\C \times \C \times \C) \]
with the factors on the left corresponding to the eigenvalues $\pm 1$ for $s$ and those on the right to the eigenvalues $\{ 1 , \rho, \rho^2 \}$ for $t$ where $\rho$ is a primitive third root of unity. The Bratelli diagrams of the inclusion of $\C$ in these factors are given by
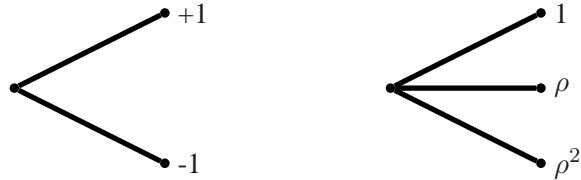
\[
\begin{pspicture}(-.5,-.5)(7.5,2.5)
\dotnode(0,1){A}
\dotnode(2,0){B}
\dotnode(2,2){C}
\dotnode(5,1){D}
\dotnode(7,0){E}
\dotnode(7,1){F}
\dotnode(7,2){G}
\uput[r](2,0){-1}
\uput[r](2,2){+1}
\uput[r](7,0){$\rho^2$}
\uput[r](7,1){$\rho$}
\uput[r](7,2){1}
\ncline[linecolor=black,linewidth=2pt]{A}{B} 
\ncline[linecolor=black,linewidth=2pt]{A}{C} 
\ncline[linecolor=black,linewidth=2pt]{D}{E} 
\ncline[linecolor=black,linewidth=2pt]{D}{F} 
\ncline[linecolor=black,linewidth=2pt]{D}{G} 
\end{pspicture}
\]
whence the structure of the quiver $\Lambda$ is given by 
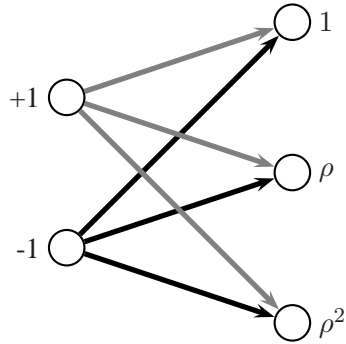
\[
\begin{pspicture}(-.5,-.5)(3.5,4.5)
\cnode(0,1){.25}{A}
\cnode(0,3){.25}{B}
\cnode(3,0){.25}{C}
\cnode(3,2){.25}{D}
\cnode(3,4){.25}{E}
\uput{10pt}[l](0,1){-1}
\uput{10pt}[l](0,3){+1}
\uput{10pt}[r](3,0){$\rho^2$}
\uput{10pt}[r](3,2){$\rho$}
\uput{10pt}[r](3,4){$1$}
\ncline[linecolor=black,linewidth=2pt]{->}{A}{C}
\ncline[linecolor=black,linewidth=2pt]{->}{A}{D}
\ncline[linecolor=black,linewidth=2pt]{->}{A}{E}
\ncline[linecolor=gray,linewidth=2pt]{->}{B}{C}
\ncline[linecolor=gray,linewidth=2pt]{->}{B}{D}
\ncline[linecolor=gray,linewidth=2pt]{->}{B}{E}
\end{pspicture}
\]
In this case the character is $\theta=(-1,-1;1,1,1)$ whence $\wis{comp}$ is the semigroup of $\Z^5$ consisting of all $(a,b;c,d,e) \in \N^5$ such that $a+b=c+d+e$ which is generated by the six dimension vectors
\[
\begin{cases}
+1, 1 = \beta_1 = (1,0;1,0,0) \\
-1 ,\rho = \beta_2 = (0,1;0,1,0) \\
+1, \rho^2 = \beta_3 = (1,0; 0,0,1) \\
-1, 1 = \beta_4 = (0,1; 1,0,0) \\
+1,\rho = \beta_5 = (1,0; 0,1,0) \\
-1, \rho^2 = \beta_6 = (0,1; 0,0,1)
\end{cases}
\]
and the quiver setting associated to $\C PSL_2(\Z)$ is $(Q,\alpha)$ with $\alpha=(1,1,1,1,1,1)$ and $Q$ the quiver
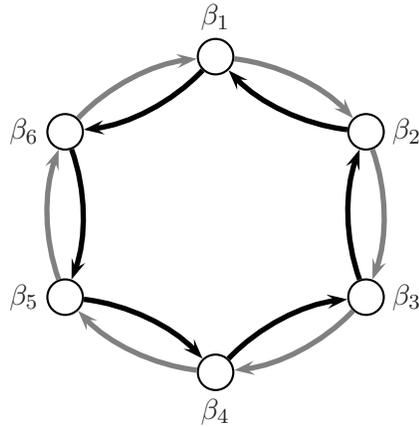
\[
\begin{pspicture}(-.5,-.5)(4.5,4.7)
\cnode(2,4.2){.25}{A}
\cnode(4,3.2){.25}{B}
\cnode(4,1){.25}{C}
\cnode(2,0){.25}{D}
\cnode(0,1){.25}{E}
\cnode(0,3.2){.25}{F}
\uput{10pt}[u](2,4.2){$\beta_1$}
\uput{10pt}[r](4,3.2){$\beta_2$}
\uput{10pt}[r](4,1){$\beta_3$}
\uput{10pt}[d](2,0){$\beta_4$}
\uput{10pt}[l](0,1){$\beta_5$}
\uput{10pt}[l](0,3.2){$\beta_6$}
\ncarc[arcangle=20,linecolor=gray,linewidth=2pt]{->}{A}{B}
\ncarc[arcangle=20,linecolor=gray,linewidth=2pt]{->}{B}{C}
\ncarc[arcangle=20,linecolor=gray,linewidth=2pt]{->}{C}{D}
\ncarc[arcangle=20,linecolor=gray,linewidth=2pt]{->}{D}{E}
\ncarc[arcangle=20,linecolor=gray,linewidth=2pt]{->}{E}{F}
\ncarc[arcangle=20,linecolor=gray,linewidth=2pt]{->}{F}{A}
\ncarc[arcangle=20,linecolor=black,linewidth=2pt]{->}{A}{F}
\ncarc[arcangle=20,linecolor=black,linewidth=2pt]{->}{F}{E}
\ncarc[arcangle=20,linecolor=black,linewidth=2pt]{->}{E}{D}
\ncarc[arcangle=20,linecolor=black,linewidth=2pt]{->}{D}{C}
\ncarc[arcangle=20,linecolor=black,linewidth=2pt]{->}{C}{B}
\ncarc[arcangle=20,linecolor=black,linewidth=2pt]{->}{B}{A}
\end{pspicture}
\]
\end{example}

As the modular group is a twofold cover of $\overline{\Gamma}$ the following example should not come as a surprise.

\begin{example}[The modular group] Because $\Gamma = \Z_4 \ast_{\Z_2} \Z_6$ we have that
\[
\C~\Gamma \simeq \C \Z_4 \ast_{\C \times \C} \C \Z_6 \]
where the left factor splits in $4$ copies of $\C$ corresponding to the eigenvalues $\{ 1,i,-1,-i \}$
and the factor on the right in $6$ copies of $\C$ corresponding to the eigenvalues $\{ 1,-\rho^2,\rho,-1,\rho^2,-\rho \}$. The amalgamating factor $\C Z_2 = \C \times \C$ has factors corresponding to the eigenvalues $\pm 1$. Hence, the Bratelli diagrams of the two inclusions are 
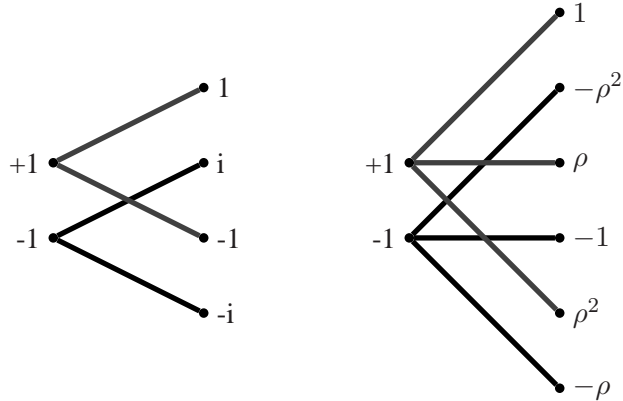
\[
\begin{pspicture}(-.5,-.5)(2.5,5.5)
\dotnode(0,2){A}
\dotnode(0,3){B}
\dotnode(2,1){C}
\dotnode(2,2){D}
\dotnode(2,3){E}
\dotnode(2,4){F}
\uput[l](0,2){-1}
\uput[l](0,3){+1}
\uput[r](2,1){-i}
\uput[r](2,2){-1}
\uput[r](2,3){i}
\uput[r](2,4){1}
\ncline[linecolor=black,linewidth=2pt]{A}{C} 
\ncline[linecolor=black,linewidth=2pt]{A}{E} 
\ncline[linecolor=darkgray,linewidth=2pt]{B}{D} 
\ncline[linecolor=darkgray,linewidth=2pt]{B}{F} 
\end{pspicture}~\qquad~\qquad
\begin{pspicture}(-.5,-.5)(2.5,5.5)
\dotnode(0,2){A}
\dotnode(0,3){B}
\dotnode(2,1){C}
\dotnode(2,2){D}
\dotnode(2,3){E}
\dotnode(2,4){F}
\dotnode(2,0){G}
\dotnode(2,5){H}
\uput[l](0,2){-1}
\uput[l](0,3){+1}
\uput[r](2,1){$\rho^2$}
\uput[r](2,2){$-1$}
\uput[r](2,3){$\rho$}
\uput[r](2,4){$-\rho^2$}
\uput[r](2,0){$-\rho$}
\uput[r](2,5){$1$}
\ncline[linecolor=black,linewidth=2pt]{A}{F} 
\ncline[linecolor=black,linewidth=2pt]{A}{D} 
\ncline[linecolor=black,linewidth=2pt]{A}{G} 
\ncline[linecolor=darkgray,linewidth=2pt]{B}{H} 
\ncline[linecolor=darkgray,linewidth=2pt]{B}{E} 
\ncline[linecolor=darkgray,linewidth=2pt]{B}{C} 
\end{pspicture}
\]
whence the quiver $\Lambda$ for $SL_2(\Z)$ is of the form
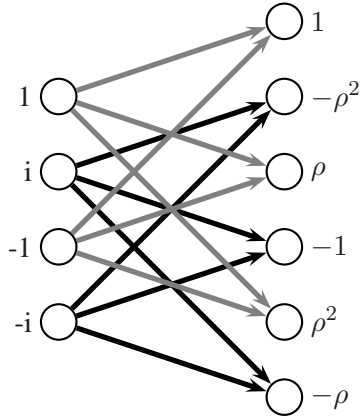
\[
\begin{pspicture}(-.5,-.5)(3.5,5.5)
\cnode(0,1){.25}{A}
\cnode(0,2){.25}{B}
\cnode(0,3){.25}{C}
\cnode(0,4){.25}{D}
\cnode(3,0){.25}{E}
\cnode(3,1){.25}{F}
\cnode(3,2){.25}{G}
\cnode(3,3){.25}{H}
\cnode(3,4){.25}{I}
\cnode(3,5){.25}{J}
\uput{10pt}[l](0,1){-i}
\uput{10pt}[l](0,2){-1}
\uput{10pt}[l](0,3){i}
\uput{10pt}[l](0,4){1}
\uput{10pt}[r](3,0){$-\rho$}
\uput{10pt}[r](3,1){$\rho^2$}
\uput{10pt}[r](3,2){$-1$}
\uput{10pt}[r](3,3){$\rho$}
\uput{10pt}[r](3,4){$-\rho^2$}
\uput{10pt}[r](3,5){$1$}
\ncline[linecolor=black,linewidth=2pt]{->}{A}{E}
\ncline[linecolor=black,linewidth=2pt]{->}{A}{G}
\ncline[linecolor=black,linewidth=2pt]{->}{A}{I}
\ncline[linecolor=black,linewidth=2pt]{->}{C}{I}
\ncline[linecolor=black,linewidth=2pt]{->}{C}{G}
\ncline[linecolor=black,linewidth=2pt]{->}{C}{E}
\ncline[linecolor=gray,linewidth=2pt]{->}{B}{F}
\ncline[linecolor=gray,linewidth=2pt]{->}{B}{H}
\ncline[linecolor=gray,linewidth=2pt]{->}{B}{J}
\ncline[linecolor=gray,linewidth=2pt]{->}{D}{F}
\ncline[linecolor=gray,linewidth=2pt]{->}{D}{H}
\ncline[linecolor=gray,linewidth=2pt]{->}{D}{J}
\end{pspicture}
\]
which we see is the disjoint union of two copies of the quiver  for $\overline{\Gamma}$. As $\theta=(-1,-1,-1,-1;1,1,1,1,1,1)$ the semigroup $\wis{comp} \subset \N^{10}$ consists of the dimension vectors $(a_1,a_2,a_3,a_4;b_1,b_2,b_3,b_4,b_5,b_6)$ such that $a_1+a_2+a_3+a_4=b_1+\hdots+b_6$ and is generated by the $12$ dimension vectors
\[
\begin{cases}
+1, 1 = \beta_1 = (1,0,0,0;1,0,0,0,0,0) \\
-1, \rho = \beta_2 = (0,0,1,0;0,0,1,0,0,0) \\
+1 ,\rho^2 = \beta_3 = (1,0,0,0; 0,0,0,0,1,0) \\
-1 ,1 = \beta_4 = (0,0,1,0; 1,0,0,0,0,0) \\
+1, \rho = \beta_5 = (1,0,0,0; 0,0,1,0,0,0) \\
-1 ,\rho^2 = \beta_6 = (0,0,1,0; 0,0,0,0,1,0)
\end{cases}~\begin{cases}
i ,-\rho^2 = \gamma_1 = (0,1,0,0;0,1,0,0,0,0) \\
-i, -1 = \gamma_2 = (0,0,0,1;0,0,0,1,0,0) \\
i ,-\rho = \gamma_3 = (0,1,0,0; 0,0,0,0,0,1) \\
-i ,-\rho^2 = \gamma_4 = (0,0,0,1; 0,1,0,0,0,0) \\
i ,-1 = \gamma_5 = (0,1,0,0; 0,0,0,1,0,0) \\
-i ,-\rho = \gamma_6 = (0,0,0,1; 0,0,0,0,0,1)
\end{cases}
\]
Hence, the quiver setting $(Q,\alpha)$ associated to $SL_2(\Z)$ is $\alpha=(1,\hdots,1)$ and $Q$ the quiver
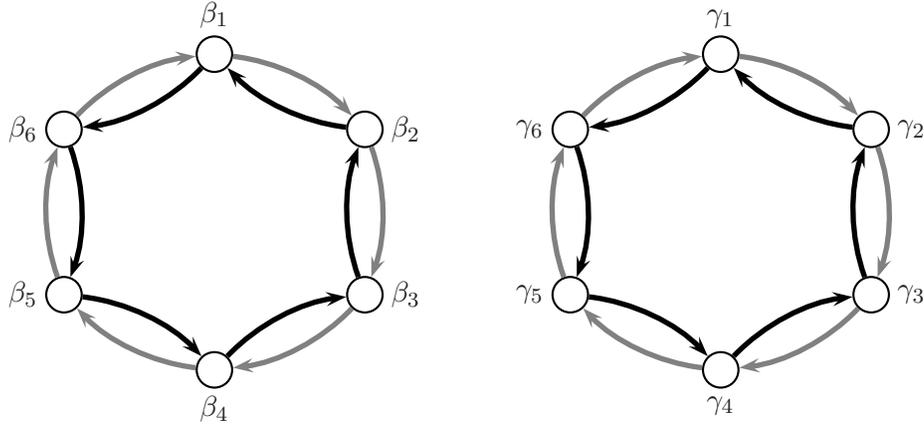
\[
\begin{pspicture}(-.5,-.5)(4.5,4.7)
\cnode(2,4.2){.25}{A}
\cnode(4,3.2){.25}{B}
\cnode(4,1){.25}{C}
\cnode(2,0){.25}{D}
\cnode(0,1){.25}{E}
\cnode(0,3.2){.25}{F}
\uput{10pt}[u](2,4.2){$\beta_1$}
\uput{10pt}[r](4,3.2){$\beta_2$}
\uput{10pt}[r](4,1){$\beta_3$}
\uput{10pt}[d](2,0){$\beta_4$}
\uput{10pt}[l](0,1){$\beta_5$}
\uput{10pt}[l](0,3.2){$\beta_6$}
\ncarc[arcangle=20,linecolor=gray,linewidth=2pt]{->}{A}{B}
\ncarc[arcangle=20,linecolor=gray,linewidth=2pt]{->}{B}{C}
\ncarc[arcangle=20,linecolor=gray,linewidth=2pt]{->}{C}{D}
\ncarc[arcangle=20,linecolor=gray,linewidth=2pt]{->}{D}{E}
\ncarc[arcangle=20,linecolor=gray,linewidth=2pt]{->}{E}{F}
\ncarc[arcangle=20,linecolor=gray,linewidth=2pt]{->}{F}{A}
\ncarc[arcangle=20,linecolor=black,linewidth=2pt]{->}{A}{F}
\ncarc[arcangle=20,linecolor=black,linewidth=2pt]{->}{F}{E}
\ncarc[arcangle=20,linecolor=black,linewidth=2pt]{->}{E}{D}
\ncarc[arcangle=20,linecolor=black,linewidth=2pt]{->}{D}{C}
\ncarc[arcangle=20,linecolor=black,linewidth=2pt]{->}{C}{B}
\ncarc[arcangle=20,linecolor=black,linewidth=2pt]{->}{B}{A}
\end{pspicture}~\qquad~\qquad
\begin{pspicture}(-.5,-.5)(4.5,4.7)
\cnode(2,4.2){.25}{A}
\cnode(4,3.2){.25}{B}
\cnode(4,1){.25}{C}
\cnode(2,0){.25}{D}
\cnode(0,1){.25}{E}
\cnode(0,3.2){.25}{F}
\uput{10pt}[u](2,4.2){$\gamma_1$}
\uput{10pt}[r](4,3.2){$\gamma_2$}
\uput{10pt}[r](4,1){$\gamma_3$}
\uput{10pt}[d](2,0){$\gamma_4$}
\uput{10pt}[l](0,1){$\gamma_5$}
\uput{10pt}[l](0,3.2){$\gamma_6$}
\ncarc[arcangle=20,linecolor=gray,linewidth=2pt]{->}{A}{B}
\ncarc[arcangle=20,linecolor=gray,linewidth=2pt]{->}{B}{C}
\ncarc[arcangle=20,linecolor=gray,linewidth=2pt]{->}{C}{D}
\ncarc[arcangle=20,linecolor=gray,linewidth=2pt]{->}{D}{E}
\ncarc[arcangle=20,linecolor=gray,linewidth=2pt]{->}{E}{F}
\ncarc[arcangle=20,linecolor=gray,linewidth=2pt]{->}{F}{A}
\ncarc[arcangle=20,linecolor=black,linewidth=2pt]{->}{A}{F}
\ncarc[arcangle=20,linecolor=black,linewidth=2pt]{->}{F}{E}
\ncarc[arcangle=20,linecolor=black,linewidth=2pt]{->}{E}{D}
\ncarc[arcangle=20,linecolor=black,linewidth=2pt]{->}{D}{C}
\ncarc[arcangle=20,linecolor=black,linewidth=2pt]{->}{C}{B}
\ncarc[arcangle=20,linecolor=black,linewidth=2pt]{->}{B}{A}
\end{pspicture}
\]
\end{example}

So far, all examples had all components of $\alpha_A$ equal to one whence we did not need a Morita equivalence. To end this section let us give an example in which the Morita equivalence does appear.

\begin{example}[Amalgamating matrices] Let $l=ak$ and $m=bk$ and consider the amalgamated coproduct
\[
A = M_l(\C) \ast_{M_k(\C)} M_m(\C) \]
The Bratelli diagrams of the two inclusions are depicted by
\[
\begin{pspicture}(-.5,-.5)(2.5,+.5)
\dotnode(0,0){A}
\dotnode(2,0){B}
\ncline[doubleline=true,linecolor=black,linewidth=2pt]{A}{B} 
\ncput*{a}
\end{pspicture}~\qquad~\qquad
\begin{pspicture}(-.5,-.5)(2.5,+.5)
\dotnode(0,0){A}
\dotnode(2,0){B}
\ncline[doubleline=true,linecolor=black,linewidth=2pt]{A}{B} 
\ncput*{b}
\end{pspicture}
\]
whence the quiver $\Lambda$ for $M_l(\C) \ast_{M_k(\C)} M_m(\C)$ is
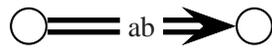
\[
\begin{pspicture}(-.5,-.5)(3.5,+.5)
\cnode(0,0){.25}{A}
\cnode(3,0){.25}{B}
\ncline[doubleline=true,linecolor=black,linewidth=2pt]{->}{A}{B} 
\ncput*{ab}
\end{pspicture}
\]
This time the character is $\theta=(-l,m)$ and the semigroup of all dimension vectors $(p,g)$ such that $lp=mq$ is generated by a unique generator
\[
\beta = (\frac{b}{(a,b)},\frac{a}{(a,b)}) \]
and as the Euler form of $\Lambda$ is given by the matrix $\begin{bmatrix} 1 & -ab \\ 0 & 1 \end{bmatrix}$ we obtain that the quiver $Q$ is
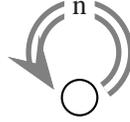
\[
\begin{pspicture}(-.5,-.5)(1,1)
\cnode(0,0){.25}{A}
\nccircle[doubleline=true,linecolor=gray,nodesep=3pt,linewidth=2pt]{->}{A}{.6cm}
\ncput*{n}
\end{pspicture}
\]
where $n = 1 + \frac{1}{(a,b)^2}(a^2b^2-a^2-b^2)$ and where $\alpha = (\frac{abk}{(a,b)})$. Hence, conjecturally, the amalgamated coproduct $M_l(\C) \ast_{M_k(\C)} M_m(\C)$ has a non-commutative \'etale extension isomorphic to
\[
M_{\frac{abk}{(a,b)}}(\C \langle x_1,\hdots,x_n \rangle) \]
\end{example}

\section{Non-commutative symplectic flows.}

The quivers constructed above are symmetric, that is, there is an involution $a \mapsto a^*$ on the arrows such that $a^*$ has the reverse orientation of $a$. Equivalently, the path algebra $\C Q$ can be equipped with a {\em symplectic structure}. In this section we will show that this is a general fact for coproducts of semi-simple algebras and illustrate how the machinery developed in \cite{LBBocklandt} can be used to construct symplectomorphisms on the quotient varieties.

\begin{proposition} Let $(Q,\alpha)$ be the quiver setting associated to the algebra coproduct
\[
S_1 \ast S_2 \]
where $S_1$ and  $S_2$ are finite dimensional semi-simple algebras. Then, the quiver $Q$ is symmetric, that is, $\C Q$ is symplectic.
\end{proposition}

\begin{proof}
If $S_1 = M_{k_1}(\C) \oplus \hdots \oplus M_{k_r}(\C)$ and $S_2 = M_{l_1}(\C) \oplus \hdots \oplus M_{l_s}(\C)$, then the Euler-form of the bipartite quiver $\Lambda$ is determined by the matrix
\[
\begin{bmatrix}
\begin{array}{ccc|ccc}
1 & & & -k_1l_1 & \hdots & -k_1l_s \\
& \ddots & & \vdots & & \vdots \\
& & 1 & -k_rl_1 & \hdots & -k_rl_s \\
\hline
0 & \hdots & 0 & 1 & & \\
\vdots & & \vdots & & \ddots & \\
0 & \hdots & 0 & & & 1
\end{array}
\end{bmatrix}
\]
Now, if $\beta=(a_1,\hdots,a_r,b_1,\hdots,b_s)$ and $\beta'=(c_1,\hdots,c_r,d_1,\hdots,d_s)$ are two generators of the semigroup $\wis{comp}$, then one computes that
\[
\chi_Q(\beta,\beta') = \sum_{i=1}^r a_ic_i + \sum_{j=1}^t b_jd_j - \sum_{i=1}^r \sum_{j=1}^s k_il_j a_id_j \]
whence
\[
\chi_Q(\beta,\beta') - \chi_Q(\beta',\beta) = \sum_{i=1}^r \sum_{j=1}^s k_il_j(a_id_j-c_ib_j) \]
Now, the components of $\beta$ and $\beta'$ must satisfy the numerical restrictions
\[
\sum_{i=1}^r k_ia_i = \sum_{j=1}^s l_j b_j \qquad \text{and} \qquad
\sum_{j=1}^r k_ic_i = \sum_{j=1}^s l_j d_j \]
from which it follows that the above expression is zero, whence there are the same number of arrows in both directions between any two distinct vertices.
\end{proof}

Observe that we only claim that the number of arrows between two different vertices is equal in each direction. It is, in general, not true that the number of loops in a vertex is even. For example, using the last example of the foregoing section one verifies that the quiver corresponding to $M_2(\C) \ast M_2(\C)$ is the three loop quiver.

We recall some facts on non-commutative differential forms for path algebras of quivers and refer for full details to \cite{LBBocklandt}. If $Q$ is a quiver on $k$ vertices, then a basis for the relative (with respect to the vertex-algebra $\C^{\times k}$) differential $n$-forms $\Omega^n_{v}~Q$ is given by the elements
\[
p_0dp_1\hdots dp_n \]
where the $p_i$ are oriented paths in $Q$ of length $\geq 1$ for $i \geq 1$ and of length $\geq 0$ for $p_0$ and such that the starting point of $p_i$ is the end point of $p_{i+1}$. Dividing out super-commutators one obtains the relative Karoubi complex
\[
\wis{dR}^0_{v}~Q \rTo^d \wis{dR}^1_{v}~Q \rTo^d \wis{dR}^2_{v}~Q \rTo^d \hdots \]
where
\[
\wis{dR}^n_{v}~Q = \frac{\Omega^n_{v}~Q}{\sum_{i=1}^n [ \Omega^i_{v}~Q,\Omega^{n-i}_{v}~Q]} \]
and one can make the first terms of this sequence quite explicit. $\wis{dR}^0_{v}~Q$ is the vectorspace spanned by all {\em necklace words} in $Q$, that is, all oriented circuits in $Q$ taken upto cyclic equivalence. $\wis{dR}^1_{v}~Q$ is the vectorspace
\[
\oplus_a p da \]
where the sum is taken over all arrows $a$ in $Q$ and where $p$ is an oriented path in $Q$ starting in the terminal vertex of $a$. The differential $d$ takes a necklace $n$ in $Q$ and takes it to 
\[
dn = \sum_a \frac{\partial n}{\partial a} da \]
where $\frac{\partial n}{\partial a}$ is the sum of all oriented paths (starting in the end point of $a$) which are obtained by deleting one occurrence  of $a$ in $n$. 

As the relative Karoubi complex is exact everywhere but the first term, we get an exact sequence
\[
0 \rTo \C^{\times k} \rTo \wis{dR}^0_{v}~Q \rTo^d \wis{dR}^1_{v}~Q_{ex} \rTo 0 \]
see \cite[\S 3]{LBBocklandt} where the last term is the kernel of $d~:~\wis{dR}^1_v~Q \rTo \wis{dR}^2_v~Q$.

In case the quiver $Q$ is {\em symplectic} (that is if we have an involution $a \mapsto a^*$ on the arrows of $Q$)  the situation is even nicer. Take a representant $a$ in each equivalence class (the other element will be denoted by $a^*$) and call this set of arrows $L$. Then, one can define a Lie algebra structure on the space $\wis{dR}^0_{v}~Q$ of non-commutative functions by
\[
[n,n'] = \sum_{a \in L} (\frac{\partial n}{\partial a}\frac{\partial n'}{\partial a^*} - \frac{\partial n}{\partial a^*}\frac{\partial n'}{\partial a}) \]
the so called {\em necklace Lie algebra}, see \cite{LBBocklandt}. Moreover, in this case $\wis{dR}^1_{v}~Q_{ex}$ can be identified with $\wis{Der}_{\omega}~Q$, the Lie algebra of {\em symplectic derivations} of $\C Q$. That is, $\C^{\times k}$-derivations of $\C Q$ which preserve the so called {\em moment-element}
\[
m = \sum_{a \in L}[a,a^*] \in \C Q \]
The identification assigns to the image $dn$ of a necklace $n$ of $Q$ the $\C^{\times k}$-derivation $\theta_n$ defined by
\[
\begin{cases}
\theta_n(a) = - \frac{\partial n}{\partial a^*} \\
\theta_n(a^*) = \frac{\partial n}{\partial a}
\end{cases}
\]
Under these identifications, the following sequence is an exact sequence of Lie algebras if we put on $\C^{\times k}$ the Abelian Lie algebra structure, see \cite[Thm. 4.2]{LBBocklandt}
\[
0 \rTo \C^{\times k} \rTo \wis{dR}^0_v~Q \rTo \wis{Der}_{\omega}~Q \rTo 0 \]
If $\theta$ is a nilpotent symplectic derivation, it determines an automorphism of the path algebra $\C Q$ and hence, by taking traces, a symplectomorphism on the quotient varieties $\wis{iss}_{\alpha}~Q$. A particular example of such derivations correspond to $dn$ where $n$ is a necklace involving only arrows from $L$.

We will give the physical relevant example of {\em Calogero-Moser particles} in which this procedure leads to the {\em Calogero-Moser flows} of \cite{Wilson}. Next, we will compute some symplectomorphisms for the quotient varieties of semi-simple representations of the modular group.

\begin{example}[Calogero-Moser flows] Consider the symplectic quiver
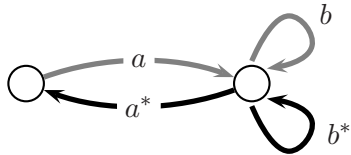
\[
\begin{pspicture}(0,0)(4,2)
\cnode(0,1){.25}{A}
\cnode(3,1){.25}{B}
\ncarc[arcangle=20,linecolor=gray,linewidth=2pt]{->}{A}{B}
\ncput*{$a$}
\ncarc[arcangle=20,linecolor=black,linewidth=2pt]{->}{B}{A}
\ncput*{$a^*$}
\pscurve[linecolor=gray,linewidth=2pt]{->}(3,1.3)(3.5,1.9)(3.8,1.5)(3.2,1.2)
\put(3.9,1.8){$b$}
\pscurve[linecolor=black,linewidth=2pt]{->}(3,0.7)(3.5,0.1)(3.8,0.5)(3.2,0.8)
\put(4,0.2){$b^*$}
\end{pspicture}
\]
The only necklaces in $Q$ involving only arrows from $L = \{ a,b \}$ are polynomials in $b$. The symplectic derivation $\theta_n$ corresponding to $n = \frac{1}{n+1} b^{n+1}$ is defined by
\[
\begin{cases}
\theta_n(a) = 0 = \theta_n(a^*) \\
\theta_n(b) = 0 \\
\theta_n(b^*) = b^n
\end{cases}
\]
and is therefore nilpotent. The corresponding automorphism $\gamma_n = e^{t \theta_n}$ is defined by
\[
\gamma_n(a)=a \quad \gamma_n(a^*)=a^* \quad \gamma_n(b) = b \quad \text{and} \quad \gamma_n(b^*) = b^* + t b^n \]
In this case, the moment-element $m=[a,a^*]+[b,b^*] = -a^*a e_1 + (aa^*+[b,b^*])e_2$. Take as the dimension vector $\alpha = (1,k)$ and consider the {\em deformed preprojective algebra}
\[
\Pi_{\lambda} = \frac{\C Q}{(m-\lambda)} \]
where $\lambda=-ke_1+1e_2$, then $\wis{iss}_{\alpha}~\Pi_{\lambda}$, the quotient variety of $\wis{rep}_{\alpha}~\Pi_{\lambda}$ under the base-change group $GL(\alpha)$, is Calogero-Moser phase space $Calo_k$. It can be identified with pairs of $k \times k$ matrices $(X,Z) \in M_k(\C) \times M_k(\C)$ (corresponding to the arrows $b^*$ resp. $b$ in $Q$) satisfying
\[
[X,Z] + 1_k \quad \text{has rank} \quad 1 \]
see \cite{Wilson} or \cite{LBBocklandt} for more details. As the coordinate ring of the quotient variety $\wis{iss}_{\alpha}~Q$ (and hence of the quotient $\wis{iss}_{\alpha}~\Pi_{\lambda}$ is generated by traces along oriented cycles in $Q$ we see that the the automorphism $\gamma_n$ defines a $GL_n$-invariant flow on $\wis{rep}_{\alpha}~\Pi_{\lambda}$ defined by
\[
(X,Z,u,v) \mapsto (X+tZ^n,Z,u,v) \]
which generate the {\em Calogero-Moser flows} on $Calo_k$ as defined in \cite{Wilson}.
\end{example}

Motivated by this example, we will now investigate some symplectic flows for the (projective) modular group $\overline{\Gamma}$, the case of $\Gamma$ being similar.

\begin{example}[The projective modular group] Let us define a symplectic structure on the quiver $Q$ corresponding to $\C \overline{\Gamma}$
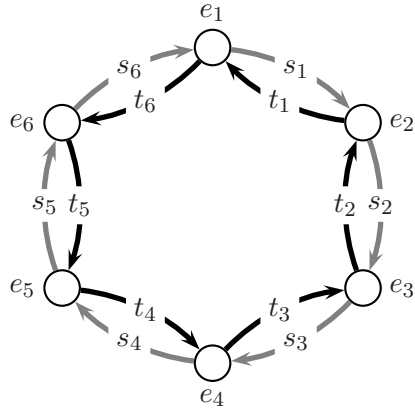
\[
\begin{pspicture}(-.5,-.5)(4.5,4.7)
\cnode(2,4.2){.25}{A}
\cnode(4,3.2){.25}{B}
\cnode(4,1){.25}{C}
\cnode(2,0){.25}{D}
\cnode(0,1){.25}{E}
\cnode(0,3.2){.25}{F}
\uput{10pt}[u](2,4.2){$e_1$}
\uput{10pt}[r](4,3.2){$e_2$}
\uput{10pt}[r](4,1){$e_3$}
\uput{10pt}[d](2,0){$e_4$}
\uput{10pt}[l](0,1){$e_5$}
\uput{10pt}[l](0,3.2){$e_6$}
\ncarc[arcangle=20,linecolor=gray,linewidth=2pt]{->}{A}{B}
\ncput*{$s_1$}
\ncarc[arcangle=20,linecolor=gray,linewidth=2pt]{->}{B}{C}
\ncput*{$s_2$}
\ncarc[arcangle=20,linecolor=gray,linewidth=2pt]{->}{C}{D}
\ncput*{$s_3$}
\ncarc[arcangle=20,linecolor=gray,linewidth=2pt]{->}{D}{E}
\ncput*{$s_4$}
\ncarc[arcangle=20,linecolor=gray,linewidth=2pt]{->}{E}{F}
\ncput*{$s_5$}
\ncarc[arcangle=20,linecolor=gray,linewidth=2pt]{->}{F}{A}
\ncput*{$s_6$}
\ncarc[arcangle=20,linecolor=black,linewidth=2pt]{->}{A}{F}
\ncput*{$t_6$}
\ncarc[arcangle=20,linecolor=black,linewidth=2pt]{->}{F}{E}
\ncput*{$t_5$}
\ncarc[arcangle=20,linecolor=black,linewidth=2pt]{->}{E}{D}
\ncput*{$t_4$}
\ncarc[arcangle=20,linecolor=black,linewidth=2pt]{->}{D}{C}
\ncput*{$t_3$}
\ncarc[arcangle=20,linecolor=black,linewidth=2pt]{->}{C}{B}
\ncput*{$t_2$}
\ncarc[arcangle=20,linecolor=black,linewidth=2pt]{->}{B}{A}
\ncput*{$t_1$}
\end{pspicture}
\]
by taking $L = \{ s_1,\hdots,s_6 \}$ and $s_i^* = t_i$. Under composition the necklaces of $Q$ are generated by the elementary necklaces
\[
C_i = s_it_i \qquad S = s_1s_2s_3s_4s_5s_6 \qquad T = t_1t_2t_3t_4t_5t_6 \]
The Lie brackets between these elementary necklaces are
\[
[C_i,C_j] = 0 \qquad [C_i,S] = -S \qquad [C_i,T] = T \qquad [S,S]=[S,T]=[T,T]=0 \]
The only necklaces involving only arrows from $L$ are the polynomials in $S$. The symplectic derivation corresponding to $n = \frac{1}{n+1}S^{n+1}$ is defined by
\[
\begin{cases}
\theta_n(s_i) = 0 \\
\theta_n(t_i) = s_{i-1}\hdots s_1 S^{n-1} s_6 \hdots s_{i+1}
\end{cases}
\]
whence the corresponding symplectic flows are determined by the automorphism
\[
\begin{cases}
\gamma_n(s_i) = s_i \\
\gamma_n(t_i) = t_i + t  s_{i-1}\hdots s_1 S^{n-1} s_6 \hdots s_{i+1}
\end{cases}
\]
which induce, after taking traces, symplectomorphisms on the quotient varieties $\wis{iss}_{\alpha}~Q$ describing all isomorphism classes of semi-simple $\alpha$-dimensional representations of $Q$.
\end{example}

\section{The modular group and the quotient singularity $\C^2/\Z_6$.}

If $A$ is formally smooth with associated quiver-setting $(Q_A,\alpha_A)$ then this setting contains enough information to describe the $GL_n$-\'etale structure of $\wis{rep}_n~A$ near the orbit of a semi-simple representation as well as the \'etale structure of the quotient varieties $\wis{iss}_n~A$ parametrizing isomorphism classes of semi-simple $n$-dimensional representations, see \cite{LBOneQuiver}. In particular, if $\{ \beta_1,\hdots,\beta_k \}$ are the semigroup generators of $\wis{comp}$ (determining the vertices of $Q_A$) and if $\beta = n_1 \beta_1 + \hdots + n_k \beta_k$ is a dimension vector for $A$, then the connected component $\wis{rep}_{\beta}~A$ contains (an open set of) simple representations if and only if $(n_1,\hdots,n_k)$ is the dimension vector of a simple representation of $Q_A$ (and there is an characterization of such dimensions vectors by \cite{LBProcesi}). For an application of this to the representation theory of $\overline{\Gamma}$ (or $B_3$) we refer to \cite{AdriLB}. Here we state the result for the modular group  $\Gamma$ :

\begin{proposition}
Let $\wis{rep}_{\alpha}~\C \Gamma$ be a connected component of $\wis{rep}_n~\C \Gamma$ such that 
\[
\alpha = b_1 \beta_1 + \hdots + b_6 \beta_6 + c_1 \gamma_1 + \hdots + c_6 \gamma_6 \]
Then, $\wis{rep}_{\alpha}~\C \Gamma$ contains simple representations if and only if one of the following two cases occurs
\[
\begin{cases}
\text{all $c_i = 0$ and $b_i \leq b_{i-1}+b_{i+1}$ for all $i~mod~6$} \\
\text{all $b_i = 0$ and $c_i \leq c_{i-1}+c_{i+1}$  for all $i~mod~6$}
\end{cases}
\]
\end{proposition}

Rather than stating the \'etale local structure results in full generality (for which we refer to \cite{AdriLB} or to \cite{LBnagatn}), we will illustrate them in some examples.

\begin{example}[Infinite dihedral group] Let $D_{\infty} = \langle s,t~|~s^2=1=t^2 \rangle$. There are four isomorphism classes of simple one-dimensional $D_{\infty}$-representations : $I_{++},I_{+-},I_{-+}$ and $I_{--}$ with for example
$I_{-+} = \C v$ where $s.v=-v$ and $t.v=v$.  As the dimension vectors $\beta_1,\hdots,\beta_4$ (corresponding to these one-dimensional simples) generated $\wis{comp}~D_{\infty}$ we know that any non-empty connected component $\wis{rep}_{\beta}~\C D_{\infty}$ of $\wis{rep}_n~\C D_{\infty}$ contains a semi-simple representation with decomposition
\[
M = I_{++}^{\oplus a_1} \oplus I_{--}^{\oplus a_2} \oplus I_{+-}^{\oplus b_1} \oplus I_{-+}^{\oplus b_2} \]
with clearly $a_1+a_2+b_1+b_2 = n$ (observe that as the $\beta_i$ are linearly dependent, there are often several descriptions $\beta = a_1 \beta_1 + a_2 \beta_4 + b_1 \beta_2 + b_2 \beta_3$). If $(Q,\alpha)$ is the quiver setting associated to $\C D_{\infty}$, then $\gamma = (a_1,a_2,b_1,b_2)$ is a dimension vector for $Q$. 

The relevance of the quiver-representation space $\wis{rep}_{\gamma}~Q$ is that it coincides as a $Stab(M) = GL(\gamma) = GL_{a_1} \times GL_{a_2} \times GL_{b_1} \times GL_{b_2}$-representation with the normal space to the orbit 
$GL_n.M = \Oscr(M)$ (which identifies as $Ext^1_{\C D_{\infty}}(M,M)$). By the \'etale slice theorems, see \cite{Luna}, we deduce that there is a $GL_n$-\'etale isomorphism between
\[
\begin{cases}
\text{$\wis{rep}_{\beta}~\C D_{\infty}$ in a neighborhood of $\Oscr(M)$, and} \\
\text{$GL_n \times^{GL(\gamma)} \wis{rep}_{\gamma}~Q$ in a neighborhood of $\Oscr(\overline{1_n \times 0})$} \end{cases}
\]
where $0$ is the zero representation in $\wis{rep}_{\gamma}~Q$. We can even make this map explicit. A representation $V \in \wis{rep}_{\gamma}~Q$ is given by matrices $S_1,S_2,T_1$ and $T_2$ where
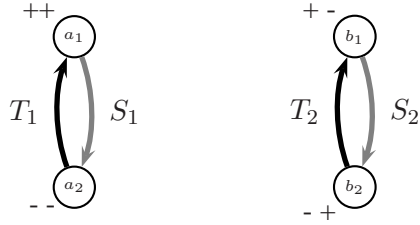
\[
\begin{pspicture}(0,-.5)(2,2.5)
\cnodeput(1,0){A}{\text{\tiny{$a_2$}}}
\uput{10pt}[dl](1,0){- -}
\cnodeput(1,2){B}{\text{\tiny{$a_1$}}}
\uput{10pt}[ul](1,2){++}
\ncarc[arcangle=20,linecolor=gray,linewidth=2pt]{->}{B}{A}
\naput{$S_1$}
\ncarc[arcangle=20,linecolor=black,linewidth=2pt]{->}{A}{B}
\naput{$T_1$}
\end{pspicture}~\qquad~\qquad
\begin{pspicture}(0,-.5)(2,2.5)
\cnodeput(1,0){A}{\text{\tiny{$b_2$}}}
\uput{10pt}[dl](1,0){- +}
\cnodeput(1,2){B}{\text{\tiny{$b_1$}}}
\uput{10pt}[ul](1,2){+ -}
\ncarc[arcangle=20,linecolor=gray,linewidth=2pt]{->}{B}{A}
\naput{$S_2$}
\ncarc[arcangle=20,linecolor=black,linewidth=2pt]{->}{A}{B}
\naput{$T_2$}
\end{pspicture}
\]
The map $\wis{rep}_{\gamma}~Q \rTo \wis{rep}_{\beta}~\C D_{\infty}$ is than given by sending the representation $V=(S_,T_1,S_2,T_2)$ to the representation determined by the matrices
\[
s \mapsto \begin{bmatrix}
\begin{array}{c c | c c}
1_{a_1} & S_1 & 0 & 0 \\
0 & -1_{a_2} & 0 & 0 \\
\hline
0 & 0 & 1_{b_1} & S_2 \\
0 & 0 & 0 & -1_{b_2}
\end{array}
\end{bmatrix}~\qquad
t \mapsto  \begin{bmatrix}
\begin{array}{c c | c c}
1_{a_1} & 0 & 0 & 0 \\
T_1 & -1_{a_2} & 0 & 0 \\
\hline
0 & 0 & -1_{b_1} & 0 \\
0 & 0 & T_2 & 1_{b_2}
\end{array}
\end{bmatrix}
\]
Therefore, in a neighborhood of $M$ any $n$-dimensional representation of $D_{\infty}$ can be conjugated to one in this matrixform.

Rather than studying these maps one at a time, we would like them to be induced from a map on the level of the algebras. The matrix-form suggest the map
\[
\C D_{\infty} \rTo^{\phi} \C Q \qquad \text{defined by} \qquad \begin{cases}
\phi(s) = e_1 -e_2 + e_3 - e_4 + s_1 + s_2 \\
\phi(t) = e_1-e_2-e_3+e_4 + t_1 +t_2
\end{cases}
\]
where $e_i$ is the vertex-idempotent of the path algebra $\C Q$ corresponding to the $i$-th vertex of $Q$ and where $s_i$ resp. $t_j$ are the arrow-variables corresponding to the red (resp. blue) arrows in $Q$. One verifies that $\phi$ is an algebra morphism.
\end{example}

In the case of the modular group we can also define an algebra map $\C \Gamma \rTo \C Q$ giving a covering of $\wis{iss}_{\alpha}~\C \Gamma$ near certain semi-simple representations. We will give the relevant information for the projective modular group $\overline{\Gamma}$, the case for $\Gamma$ follows by taking block matrices as in the case of $D_{\infty}$.

\begin{example}[The projective modular group] Let us take the notation for vertices and arrows as in the foregoing section. 
Fixing a dimension vector $\gamma=(a_1,\hdots,a_6)$ we can write the \'etale map
\[
GL_n \times^{GL(\gamma)} \wis{rep}_{\gamma}~Q \rTo \wis{rep}_n~\C \overline{\Gamma} \]
in a neighborhood of the semi-simple $n$-dimensional representation 
\[
I_{1,1}^{\oplus a_1} \oplus \hdots \oplus I_{-1,\rho^2}^{\oplus a_6}
\]
 as induced by the map $\wis{rep}_{\gamma}~Q \rTo \wis{rep}_n~\C \overline{\Gamma}$ sending a representation $V=(S_1,\hdots,S_6,T_1,\hdots,T_6)$ (with obvious notation) to the representation determined by the $n \times n$ matrices
\[
s \mapsto \begin{bmatrix}
1_{a_1} & i \rho S_1 & 0 & 0 & 0 & T_6 \\
0 & -1_{a_2} & 0 & 0 & 0 & 0 \\
0 & T_2 & 1_{a_3} & i \rho^2 S_3 & 0 & 0 \\
0 & 0 & 0 & -1_{a_4} & 0 & 0 \\
0 & 0 & 0 & T_4 & 1_{a_5} & iS_5 \\
0 & 0 & 0 & 0 & 0 & -1_{a_6}
\end{bmatrix}
\]
\[
t \mapsto \begin{bmatrix}
1_{a_1} & 0 & 0 & 0 & 0 & 0 \\
iT_1 & \rho 1_{a_2} & S_2 & 0 & 0 & 0 \\
0 & 0 & \rho^2 1_{a_3} & 0 & 0 & 0 \\
0 & 0 & iT_3 & 1_{a_4} & \rho S_4 & 0 \\
0 & 0 & 0 & 0 & \rho 1_{a_5} & 0 \\
\rho^2 S_6 & 0 & 0 & 0 & iT_5 & \rho^2 1_{a_6}
\end{bmatrix}
\]
These matrices suggest the following map $\phi~:~\C \overline{\Gamma} \rTo \C Q$ given by
\[
\begin{cases}
\phi(s) = e_1-e_2+e_3-e_4+e_5-e_6+i \rho s_1 + i\rho^2 s_3 + i s_5 + t_2+t_4+t_6 \\
\phi(t) = e_1 + \rho e_2 + \rho^2 e_3 + e_4 + \rho e_5 + \rho^2 e_6 + s_2 + \rho s_4 + \rho^2 s_6 + it_1 + it_3+it_5
\end{cases}
\]
which one verifies to be an algebra map. The specific choice of the scalar factors has the advantage that we get slightly nicer matrices for $\sigma_1 =  t^{-1}s = t^2 s$, that is, $V(\sigma_1)=$
\[
\begin{bmatrix}
1_{a_1} & i \rho S_1 & 0 & 0 & 0 & T_6 \\
-i\rho^2 T_1 & S_1T_1-T_2 S_2-\rho^2 1_{a_2} & -S_2 & -i\rho^2 S_3S_2 & 0 & -i \rho^2 T_6T_1 \\
0 & \rho T_2 & \rho 1_{a_3} & i S_3 & 0 & 0 \\
0 & -i\rho T_2T_3 & -i\rho T_3 & S_3T_3-T_4S_4 -1_{a_4} & -S_4 & -i S_5S_4 \\
0 & 0 & 0 & \rho^2 T_4 & \rho^2 1_{a_5} & i \rho^2 S_5 \\
-S_6 & -i \rho S_1S_6 & 0 & -i  T_4T_5 & -i T_5 & S_5T_5-T_6S_6-\rho 1_{a_6}
\end{bmatrix}
\]
and $\sigma_2 = st^{-1} = st^2$ whence $V(\sigma_2) =$
\[
\begin{bmatrix}
T_1S_1-S_6T_6+1_{a_1} & i S_1 & -i \rho S_2S_1 & 0 & -i T_5T_6 & \rho T_6 \\
i \rho^2 T_1 & -\rho^2 1_{a_2} & S_2 & 0 & 0 & 0 \\
-i \rho^2 T_1T_2 & \rho^2 T_2 & T_3S_3-S_2T_2 + \rho 1_{a_3} & i \rho^2 S_3 & -i \rho^2 S_4S_3 & 0 \\
0 & 0 & i \rho T_3 & -1_{a_4} & S_4 & 0 \\
-i S_6S_5 & 0 & -i \rho T_3T_4 & T_4 & \rho^2 1_{a_5} + T_5S_5-S_4T_4 & i \rho S_5 \\
S_6 & 0 & 0 & 0 & i T_5 & -\rho 1_{a_6} 
\end{bmatrix}
\]
Summarizing our arguments and using the quotient map $B_3 \rOnto B_3/\langle c \rangle = \overline{\Gamma}$ we obtain the following result on finite dimensional representations of the third braid group.

\begin{proposition} Let $\gamma=(a_1,\hdots,a_6)$ be a dimension vector for $Q$ and let 
\[
V=(S_1,\hdots,S_6,T_1,\hdots,T_6) \in \wis{rep}_{\gamma}~Q \]
Then, for any scalar $\lambda \in \C^*$ we have that
\[
\sigma_1 \mapsto \lambda V(\sigma_1) \qquad \text{and} \qquad \sigma_2 \mapsto \lambda V(\sigma_2) \]
(with the matrices $V(\sigma_i)$ given above) is an $n=\sum_i a_i$-dimensional representation of the third braid group $B_3$.

Moreover, if for all $i~mod~6$ we have that $a_i \leq a_{i-1}+a_{i+1}$ then this representation is irreducible for general $V$ and finally, a sufficiently general irreducible representation of $B_3$ can be conjugated into such form for suitable dimensions $a_1,\hdots,a_6$ satisfying $a_i \leq a_{i-1}+a_{i+1}$ for all $i~mod~6$.
\end{proposition}

Two simple $B_3$-representations with the same dimension vector are isomorphic if and only if their scalar factors coincide and for the representations of $Q$ all traces along necklaces are the same. Although this gives, in principle, a classification of simple $B_3$ (or $\Gamma$ or $\overline{\Gamma}$) representations (certainly in low dimensions), not much is known on the geometry of the quotient varieties $\wis{iss}_{\gamma}~Q$. For example, if $gcd(\gamma) = gcd(a_i)$ is not a divisor of 420, it is not known whether this variety is stably rational.

For this reason it is important to construct families of classifiable representations of $Q$ (and hence of the groups mentioned). We will show how the symplectic structure on the quiver $Q$ allows us to find such families corresponding to points in the quotient singularity $\C^2/\Z_6$ and its blow-ups.

For the symplectic structure on the quiver $Q$ given above, the corresponding moment-element
$
m = \sum_{i=1}^6 [s_i,t_i] =$
\[
 (s_6t_6-t_1s_1)e_1 + (s_1t_1-t_2s_2)e_2+(s_2t_2-t_3e_3)e_3+\hdots+(s_5t_5-t_6e_6)e_6 \]
which explains the strange looking scalar factors in the matrix descriptions of the images of $s$ and $t$ given before as for these we encounter these terms on the main diagonal for the corresponding images of $\sigma_1$ and $\sigma_2$. Consider the following two central elements of $\C Q$
\[
\begin{cases}
\lambda_0 = 0  \\
\lambda_1 = 1 e_1 + \rho^2 e_2 + \rho e_3 + 1 e_4 + \rho^2 e_5 + \rho e_6
\end{cases}
\]
and the corresponding (deformed) preprojective algebras
\[
\Pi_0 = \frac{\C Q}{(m-\lambda_0)} \qquad \text{and} \qquad \Pi_1 = \frac{\C Q}{(m- \lambda_1)}
\]
First, consider the dimension vector (the isotropic Schur root of the tame quiver $\tilde{A}_5$) $\delta = (1,1,1,1,1,1)$ and observe that $\lambda_1.\delta = 0$. It is well known that
\[
\wis{iss}_{\gamma_1}~\Pi_0 = \C^2/\Z_6 \]
where the isolated singularity corresponds to the zero-representation of $Q$ and that
\[
\wis{iss}_{\delta}~\Pi_1 = N(\lambda_1,\delta) \]
is a two-dimensional affine variety which is a deformation of $\C^2/\Z_6$ (the so-called Marsden-Weinstein reduction), see for example \cite{Kronheimer}.

For an arbitrary dimension vector $\gamma = (a_1,\hdots,a_6)$ we can describe the quotient varieties $\wis{iss}_{\gamma}~\Pi_0$ and (when applicable) $\wis{iss}_{\gamma}~\Pi_1$ using the results of \cite{Crawley}. If $a = max(a_1,\hdots,a_6)$, then
\[
\wis{iss}_{\gamma}~\Pi_0 \simeq S^a(\C^2/\Z_6) \]
the $a$-the symmetric power of the Kleinian singularity $\C^2/\Z_6$ embedded in $\wis{iss}_{\gamma}~Q$. The closed subvariety of $\wis{rep}_{\alpha}~\overline{\Gamma}$ corresponding to the image of $\wis{rep}_{\gamma}~\Pi_0$ under $\phi^*$ consists of matrix-pairs
\[ \sigma_1 \mapsto 
\begin{bmatrix}
1_{a_1} & i \rho S_1 & 0 & 0 & 0 & T_6 \\
-i\rho^2 T_1 & -\rho^2 1_{a_2} & -S_2 & -i\rho^2 S_3S_2 & 0 & -i \rho^2 T_6T_1 \\
0 & \rho T_2 & \rho 1_{a_3} & i S_3 & 0 & 0 \\
0 & -i\rho T_2T_3 & -i\rho T_3 &  -1_{a_4} & -S_4 & -i S_5S_4 \\
0 & 0 & 0 & \rho^2 T_4 & \rho^2 1_{a_5} & i \rho^2 S_5 \\
-S_6 & -i \rho S_1S_6 & 0 & -i  T_4T_5 & -i T_5 & -\rho 1_{a_6}
\end{bmatrix}
\]
\[ \sigma_2 \mapsto
\begin{bmatrix}
1_{a_1} & i S_1 & -i \rho S_2S_1 & 0 & -i T_5T_6 & \rho T_6 \\
i \rho^2 T_1 & -\rho^2 1_{a_2} & S_2 & 0 & 0 & 0 \\
-i \rho^2 T_1T_2 & \rho^2 T_2 &  \rho 1_{a_3} & i \rho^2 S_3 & -i \rho^2 S_4S_3 & 0 \\
0 & 0 & i \rho T_3 & -1_{a_4} & S_4 & 0 \\
-i S_6S_5 & 0 & -i \rho T_3T_4 & T_4 & \rho^2 1_{a_5} & i \rho S_5 \\
S_6 & 0 & 0 & 0 & i T_5 & -\rho 1_{a_6} 
\end{bmatrix}
\]
with $(S_1,\hdots,S_6,T_1,\hdots,T_6) \in \wis{rep}_{\gamma}~\Pi_0$. For special dimension vectors $\gamma$ (those satisfying $\lambda_1.\gamma = 0$), there will also be a component $\wis{iss}_{\gamma}~\Pi_1$ which is isomorphic to
\[
\wis{iss}_{\gamma}~\Pi_1 \simeq S^a(N(\lambda_1,\delta)) \]
and the closed subvariety of $\wis{iss}_{\alpha}~\overline{\Gamma}$ corresponding to the image of $\wis{rep}_{\gamma}~\Pi_1$ under $\phi^*$ consists of matrix-pairs
 \[
 \sigma_1 \mapsto \begin{bmatrix}
1_{a_1} & i \rho S_1 & 0 & 0 & 0 & T_6 \\
-i\rho^2 T_1 & 0 & -S_2 & -i\rho^2 S_3S_2 & 0 & -i \rho^2 T_6T_1 \\
0 & \rho T_2 & \rho 1_{a_3} & i S_3 & 0 & 0 \\
0 & -i\rho T_2T_3 & -i\rho T_3 & 0  & -S_4 & -i S_5S_4 \\
0 & 0 & 0 & \rho^2 T_4 & \rho^2 1_{a_5} & i \rho^2 S_5 \\
-S_6 & -i \rho S_1S_6 & 0 & -i  T_4T_5 & -i T_5 & 0
\end{bmatrix}
\]
\[
\sigma_2 \mapsto
\begin{bmatrix}
0 & i S_1 & -i \rho S_2S_1 & 0 & -i T_5T_6 & \rho T_6 \\
i \rho^2 T_1 & -\rho^2 1_{a_2} & S_2 & 0 & 0 & 0 \\
-i \rho^2 T_1T_2 & \rho^2 T_2 & 0 & i \rho^2 S_3 & -i \rho^2 S_4S_3 & 0 \\
0 & 0 & i \rho T_3 & -1_{a_4} & S_4 & 0 \\
-i S_6S_5 & 0 & -i \rho T_3T_4 & T_4 & 0 & i \rho S_5 \\
S_6 & 0 & 0 & 0 & i T_5 & -\rho 1_{a_6} 
\end{bmatrix}
\]
where $(S_1,\hdots,S_6,T_1,\hdots,T_6)$ is a $6n$-dimensional representation of $\Pi_1$. 
 It would be interesting to study these families further, in particular in connection with semi-simple modular tensor categories \cite{BakalovKirillov} and the construction of a knot-invariant starting from a semi-simple $B_3$-representation as presented in \cite{Westbury2}.
\end{example}

\section{Towards non-commutative \'etale maps.}

Clearly we would like the algebra morphisms $\C D_{\infty} \rTo \C Q$ and $\C \overline{\Gamma} \rTo \C Q$ to be \'etale in some non-commutative definition. However, this definition is not entirely clear at this moment. 

The notion of non-commutative {\em formally \'etale morphism} was introduced in \cite{KontRos} extending Grothendieck's characterization of \'etale commutative morphisms. A $\C$-algebra morphism $A \rTo^{\phi} B$ is said to be formally \'etale if for every test-object $(T,I)$ (that is, $I$ is a nilpotent ideal of $T$) and every commuting diagram
\[
\begin{diagram}
T & \rOnto^{\pi} & T/I \\
\uTo^{f} & \luDotsto^{\tilde{g}} & \uTo^{g} \\
A & \rTo^{\psi} & B
\end{diagram}
\]
there exists a {\em unique} lift $\tilde{g}$ making both triangles commute. This definition, however, is too restrictive to be of any practical interest.

\begin{example} If $S = M_{n_1}(\C) \oplus \hdots \oplus M_{n_k}(\C)$ is a finite dimensional semi-simple algebra one would like to have that $\C \rInto S$ is a non-commutative \'etale morphism. As $S$ is formally smooth, there do exist lifts for every test-object
\[
\begin{diagram}
T & \rOnto^{\pi} &  T/I \\
\uTo & \luDotsto^{\tilde{g}} & \uTo^g \\
\C & \rInto & S \end{diagram}
\]
However, this lift is usually {\em not} unique as one can conjugate it with any non-central unit $t \in T^*$ such that $\pi(t)=1$. Actually, Cuntz and Quillen have proved in \cite[Prop.6.1]{CuntzQuillen} that for any two lifts $\tilde{g}$ and $\tilde{g}'$ of $g$ there exists a $t \in T^*$ such that
\[
\tilde{g}'(s) = t \tilde{g}(s) t^{-1} \]
Hence, to get a useful workable of non-commutative \'etale morphism one has to allow at least these conjugated lifts.
\end{example}

\begin{definition} An algebra map $A \rTo^{\phi} B$ is said to be a {\em non-commutative strong \'etale} morphism if there are semi-simple subalgebras $S_A \rInto A$ and $S_B \rInto B$ such that $\phi(S_A) \subset S_B$ having the following lifting property. For every test-object and every (square) commuting diagram
\[
\psmatrix
S_B & B & T/I \\
S_A & A & T
\endpsmatrix
\psset{shortput=nab,arrows=->,labelsep=3pt}
\ncline{2,1}{1,1}^{\phi}
\ncline{2,1}{2,2}
\ncline{1,1}{1,2}
\ncline{2,2}{1,2}
\naput[npos=.4]{\phi}
\ncline{1,2}{1,3}^{g}
\ncline{2,2}{2,3}
\nbput{f}
\ncline{2,3}{1,3}
\nbput{\pi}
\ncline[linecolor=gray,border=3pt]{1,1}{2,3}
\ncput*[npos=.6]{\tilde{f_s}}
\ncline[linecolor=darkgray]{1,2}{2,3}
\naput{\exists !}
\]
and every lift $\tilde{f_s}~:~S_B \rTo T$ making the triangles $(S_A,S_B,T)$ and $(S_B,T,T/I)$ commute there is a {\em unique} lift $\tilde{f}~:~B \rTo T$ such that the triangles $(A,B,T)$ and $(B,T,T/I)$ are commutative.
\end{definition}

Whereas this definition accommodates for the problem mentioned in the example, we believe it is {\em not} the definite definition of a {\em non-commutative \'etale morphism} but it will serve our present purposes. Indeed, the algebra maps
\[
\C D_{\infty} \rTo \C Q \qquad \text{and} \qquad \C \overline{\Gamma} \rTo \C Q \]
constructed above are strong \'etale taking for the semi-simple algebras $\C$ (in $\C D_{\infty}$ and $\C \overline{\Gamma}$) and the vertex-semisimple algebra in $\C Q$. The uniqueness of the lift is then imposed by the fact that there is at most one arrow in a given direction between any two vertices in $Q$. As a consequence, the image of the arrow-variable can be obtained from the images of $s$ and $t$ and the vertex-idempotents $e_i$.

Perhaps the definite theory of non-commutative \'etale morphisms will be a non-commutative Grothendieck topology in the definition of Fred Van Oystaeyen \cite{FVOAGAA}. We note that a truly non-commutative Zariski topology on $\wis{rep}~A$ has been defined in \cite{LBnotop}.

\end{document}